\documentclass[aps,pre,twocolumn,floats,superscriptaddress,floatfix]{revtex4}
\usepackage{amsmath,amssymb}
\usepackage{graphicx}
\usepackage{psfrag}
\usepackage{url}
\usepackage{dcolumn}
\usepackage{bm}
\renewcommand{\v}{{\bm v}}
\newcommand{\dtau}{{\partial_\tau}}

\newcommand{\gradx}{{\nabla_{\x}}}
\newcommand{\lapx}{{\nabla^2_{\x}}}

\newcommand{\phig}{\varphi_{\mathrm{g}}}
\def\ph{{\hat p}}
\def\qh{{\hat q}}
\def\rh{{\hat r}}
\def\sh{{\hat s}}
\def\th{{\hat t}}
\def\rf#1{(\ref{#1})}
\def\R{{\mathbb R}}
\def\x{\bm x}
\def\q{\bm q}
\def\du{\partial_1}
\def\dd{\partial_2}
\def\duu{\partial^2_{11}}
\def\dud{\partial^2_{12}}
\def\ddd{\partial^2_{22}}
\def\ue{\mathrm{e}}
\def\ui{\mathrm{i}}
\def\pf{\mathfrak{p}}
\begin{document}

\title{Optimal transport by omni-potential flow and cosmological reconstruction}

\author{Uriel Frisch}
\email{uriel@oca.eu}
\affiliation{UNS, CNRS, Lab. Lagrange, OCA, B.P. 4229, 06304 Nice Cedex 4, France}
\author{Olga Podvigina}
\email{olgap@mitp.ru}
\affiliation{Institute of Earthquake Prediction Theory and Mathematical Geophysics of the Russian Academy of Sciences, 84/32 Profsoyuznaya St., 117997 Moscow, Russian Federation}
\author{Barbara Villone}
\email{villone@to.infn.it}
\affiliation{INAF, Osservatorio Astrofisico di Torino, Via Osservatorio, 20, 10025 Pino Torinese, Torino, Italy}
\author{Vladislav Zheligovsky}
\email{vlad@mitp.ru}
\affiliation{Institute of Earthquake Prediction Theory and Mathematical Geophysics of the Russian Academy of Sciences, 84/32 Profsoyuznaya St., 117997 Moscow, Russian Federation}
\date{\today}

\begin{abstract}
One of the simplest models used in studying the dynamics of large-scale
structure in cosmology, known as the Zeldovich approximation, is
equivalent to the three-dimensional inviscid Burgers equation for potential
flow. For smooth initial data and sufficiently short times it has the property
that the mapping of the positions of fluid particles at any time $t_1$ to
their positions at any time $t_2\ge t_1$ is the gradient of a convex
potential, a property we call \textit{omni-potentiality}. Are there other flows
with this property, that are not straightforward generalizations of Zeldovich
flows? This is answered in the affirmative in both two and three
dimensions. How general are such flows? Using a WKB technique we show that in
two dimensions, for sufficiently short times, there are omni-potential flows with arbitrary
smooth initial velocity. Mappings with a convex potential are known to be
associated with the quadratic-cost optimal transport problem. This has important
implications for the problem of reconstructing the dynamical history of the
Universe from the knowledge of the present mass distribution.

\begin{center}
{\em Dedicated to the memory of Roman Juszkiewicz}
\end{center}
\end{abstract}

\maketitle

\section{Introduction}
\label{s:intro}

Reconstruction in cosmology considers the following problem: one assumes that
the present spatial distribution of masses (galaxies and clusters, including
their dark-matter components) is known from observations, and one wants to
reconstruct the dynamical history of the Universe all the way to the earliest
epoch, when matter and radiation decoupled (nearly 14 billion years
ago). Peebles \cite{peebles89} introduced the reconstruction problem, and
proposed a variational formulation for solving it on a relatively small spatial
scale, that of the Local Group (which includes our own galaxy and neighboring
ones). On scales much larger than that of the Local Group, which have been
mapped in recent years through various projects such as the Sloan Digital Sky
Survey \cite{sdss}, reconstruction may be posed in the simplest cases as
an instance of optimal mass transport. Indeed, Frisch et al. \cite{nature}
showed that when the \textit{Zeldovich approximation} \cite{zeldovich} or
a refinement thereof (cf.~below) are applied to the relevant cosmological fluid
equations, the correspondence between the positions of mass elements initially
(at decoupling) and finally (at the present epoch) is the solution to an optimal
mass transport problem with quadratic cost. This solution is uniquely
prescribed by the \textit{marginals}: the mass distribution at decoupling
(essentially uniform) and its highly non-uniform present
distribution. A striking feature is that the sole knowledge of the current
positions of galaxies, without knowledge of their (proper)
velocities, yields nevertheless a unique solution for this kind of
large-scale reconstruction.

It was then shown by Brenier et al. and by Loeper \cite{mnras,loeper} that,
with prescribed marginals, unique reconstruction, not only of the Lagrangian
map, but of the full dynamical history of matter elements, carries over to the
Euler--Poisson model, whose validity extends much beyond that of the Zeldovich
approximation. Its unique solution is again obtained from an optimal transport
problem with a convex cost function, expressible as a space-time integral of
a suitable action, a problem whose numerical resolution remains a challenge.

As is well known, the mass transport problem was introduced by Monge
\cite{monge1781} more than two hundred years ago, and the theory took its
modern shape after the 1942 work of Kantorovich \cite{kantorovich} (see, e.g.,
Villani \cite{villani2009} for review).

The Zeldovich approximation \cite{zeldovich} was introduced in 1970 as
a first formulation in terms of Lagrangian coordinates of the growth
of density perturbations. It replaces the full Euler--Poisson
equations by basically the three-dimensional inviscid Burgers
equation (written here in standard fluid dynamical notation)
\begin{equation}
\partial_t{\bf v}+{\bf v}\cdot\nabla{\bf v}=0, \quad {\bf v}=\nabla\varphi.
\label{3dburgers}
\end{equation}
The validity of the Zeldovich approximation is controlled by how close one is
to decoupling, but in a
scale-dependent way: at very large scales, the Zeldovich approximation remains
valid up to the present epoch; at very small scales, the formation of
multi-stream caustics quickly ruins not only the validity of the Zeldovich
approximation, but even that of the
Euler--Poisson model. An immediate consequence of \eqref{3dburgers}
is that the velocity of any fluid particle remains constant in the
course of time and that the trajectories are straight lines. We denote
by $\q$ the initial (Lagrangian) fluid particle positions and
by $\x$ their (Eulerian) positions at the current epoch $t=T$. The
Lagrangian map associated with the Zeldovich approximation is
\begin{equation}
\q\mapsto\x=\nabla_{\q}\left(\frac{|\q|^2}{2}+T\varphi_0(\q)\right),
\label{burglagmap}
\end{equation}
where $\varphi_0(\q)$ is the initial velocity potential. For sufficiently
small $T$ and a sufficiently smooth initial potential, the Lagrangian map is thus
the gradient of a convex function, a property shared by the next-order
approximation, which will be discussed in Sec.~\ref{s:cosmo}. This is why
reconstruction is linked to optimal transport; indeed,
a theorem of Brenier \cite{brenierlong} states that the solution to the Monge
optimal transport problem with quadratic cost is a gradient
of a convex function, which satisfies a Monge--Amp\`ere equation. The method
of cosmological reconstruction in which one assumes that the Lagrangian map
has a convex potential and then numerically solves a quadratic-cost optimal
transport problem (after suitable discretization) is called the
Monge--Amp\`ere--Kantorovich (MAK) method \cite{nature}.

The Zeldovich approximation gives us some insight into the full temporal history
of mass elements. An important consequence of \eqref{3dburgers} is that for any
\hbox{$0\le t_1<t_2\le T$} the mapping of particle positions at time $t_1$ to their
positions at time $t_2$ is also a gradient of a convex function. When the
flow-induced mapping between any two times is potential, the flow is here
called \textit{omni-potential}. As we shall see, the velocity field associated
with such flows has the property of being simultaneously potential in Eulerian
coordinates (in cosmology, this constraint stems from the potential character
of gravity and the expansion of the Universe), as well as in Lagrangian
coordinates (which allows reconstruction by solving an optimal transport problem).

We are of course led to ask whether there exist omni-potential flows other than
Zeldovich/Burgers ones or trivial variants thereof. Investigating this issue is
the central topic of our paper.

In Section~\ref{s:geometry} we show that omni-potentiality can be reexpressed
geometrically and algebraically in terms of Hessian matrices and recast as a
set of one or several partial differential equations (depending on the space
dimension). In Section~\ref{s:algebra} we use an algebraic method to construct
explicit non-trivial examples of omni-potential flows in two and three
dimensions. These are rather special and we are led to ask how general are
omni-potential flows. In Section~\ref{s:wkb} we construct a fairly general
class of two-dimensional omni-potential flows, leaving the problem open in
higher dimensions. In Section~\ref{s:cosmo} we return to questions of
cosmological relevance: To what extent are the full solutions to the
cosmological equations omni-potential? Why is the validity of MAK
reconstruction better than that of the Zeldovich approximation, as pointed out
by Mohayaee et al.~\cite{royaetal}? In Section~\ref{s:conclusion} we list
some open problems and make concluding remarks. Finally, in the Appendix we
characterize sets of commuting symmetric matrices by constructing
suitable invariants.

\section{Criteria for omni-potentiality of flows}
\label{s:geometry}

In the present paper we study the kinematics of omni-potential flows.
We start by recalling the basic definitions of the Lagrangian and Eulerian
description of a flow --- a~motion of fluid regarded as a continuum
of infinitesimal fluid particles (whose mathematical abstraction is
a {\em point} particle).

Denote by ${\bf v}(\x,t)$ the velocity of the fluid measured at point $\x$
in space and at time $t$. It is usually called the {\em Eulerian velocity}
(the velocity measured at a fixed position in the laboratory frame).
The motion of a fluid particle satisfies the ordinary differential equation
$$\dot\x={\bf v}(\x,t),$$
which has to be supplemented by the initial condition
$$\x|_{t=0}=\q.$$
If the velocity field ${\bf v}(\x,t)$ is prescribed and sufficiently regular,
one can solve this initial value problem, at least locally in time, and obtain
the mapping $\q\mapsto\x(\q,t)$ called the {\em Lagrangian map}.
It takes a particle at the Lagrangian position $\q$ and carries it to the
Eulerian position $\x$. For a fixed $\q$, the curve
$\x(\q,t)$, parameterized by varying time $t$, is the \textit{trajectory
of the particle}, whose Lagrangian position is $\q$.
When, in a field associated with the flow,
we perform the substitution $\x\to\x(\q,t)$, we obtain its Lagrangian
description, in which the field is now a function
of the Lagrangian coordinates $\q$. For instance,
${\bf v}(\x(\q,t),t)$ is called the {\em Lagrangian velocity}.

In this paper, we consider flows ${\bf v}(\x,t)$ defined, for simplicity, in
the entire space $\R^d$, but restricted to a finite time interval
$[0,T]$. The flow induces a set of mappings of space: given two arbitrary
times $t$ and $\tau$, such that $0\le t<\tau\le T$, the
mapping from fluid particle positions at time $t$ to their positions at time
$\tau$ is here called the $(t,\tau)$-mapping. The $(0,\tau)$-mapping is just
the standard Lagrangian map.

\begin{figure}
\includegraphics[width=85mm]{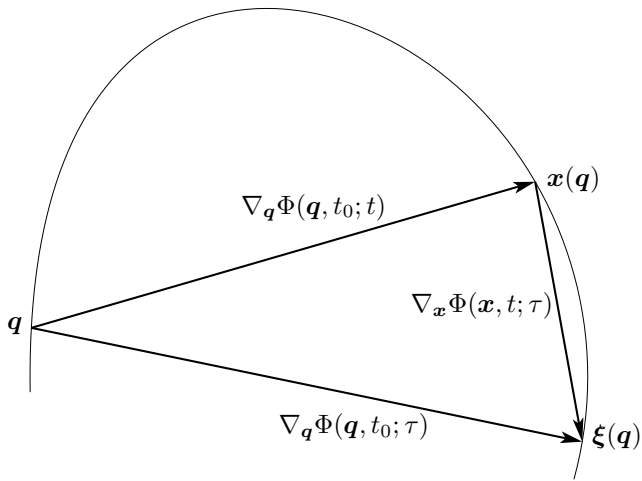}
\caption{Potential composition of two potential maps.}
\label{f:fig1}\end{figure}

As stated in the Introduction, for any two times $t$ and $\tau$, such that
$0\le t<\tau\le T$, the $(t,\tau)$-mappings induced by {\em omni-potential}
flows are required to be the gradients of a convex potential $\Phi(\q,t;\tau)$:
\begin{equation}
\q\mapsto\x=\nabla_{\q}\Phi(\q,t;\tau).
\label{potmap}\end{equation}
Such a mapping is here called potential and $\Phi(\q,t;\tau)$ is
called the $(t,\tau)$-potential. Given any three times $t_0$, $t$ and $\tau$
such that $0\le t_0\le t\le\tau\le T$, the $(t_0,t)$-mapping composed with
the $(t,\tau)$-mapping obviously yields the $(t_0,\tau)$-mapping.
This semigroup associativity, combined with omni-potentiality, implies
\begin{equation}
\Phi(\q,t_0;\tau)=\Phi\left(\nabla_{\q}\Phi(\q,t_0;t),t;\tau\right),
\label{associative}
\end{equation}
which is illustrated in Fig.~\ref{f:fig1}.

Differentiation of \eqref{potmap} with respect to time $\tau$ shows that
the velocity, ${\bf v}(\q,t;\tau)$, is also potential:
$${\bf v}(\q,t;\tau)=\nabla_{\q}\,\varphi(\q,t;\tau),$$
where
$$\varphi(\q,t;\tau)\equiv\partial_{\tau}\Phi(\q,t;\tau).$$
In this notation, ${\bf v}(\q,0;t)=\nabla_{\q}\varphi(\q,0;t)$ is
the Lagrangian velocity of a fluid particle, while its Eulerian velocity is
${\bf v}(\x,t;t)=\nabla_{\x}\varphi(\x,t;t)$; thus, in omni-potential
flow, the velocity is potential in both Eulerian and Lagrangian
coordinates. In three dimensions, this is equivalent, locally, to the
statement that the vorticity (the curl of the velocity) vanishes in
both coordinate systems. We shall call this dual potentiality in both Eulerian
and Lagrangian coordinates \textit{bi-potentiality}.

We demand that all $(t,\tau)$-potentials of omni-potential flows be convex in
the space variables. This implies that the associated $(t,\tau)$-mappings have
inverses that are also potential, the potentials associated with the inverse
mappings being the Legendre transforms of those of the direct mappings (see
Ref.~\cite{mnras}, Appendix~C and references therein). Actually,
for a sufficiently smooth flow, the converse is true: invertibility implies
convexity of the maps, as now explained. Let us denote by ${\cal H}(f)$
the Hessian of a (twice differentiable) function $f(\q)$, i.e., the $d\times d$
matrix with the entries
$${\cal H}_{ij}(f)\equiv\partial^2_{q_iq_j}f.$$
We now assume that the potentials are twice differentiable in the space
variables and that their Hessians
${\cal H}(\Phi(\q,t;\tau))$ depend continuously on their time
arguments. We then observe that, for coinciding times, the $(t,t)$-mapping is
the identity mapping, clearly having the convex potential
$$\Phi(\q,t;t)=|\q|^2/2,$$
whose Hessian is the identity. As the two times separate, loss of convexity
would require one or several eigenvalues of the Hessian to change sign and thus
to go through zero; at such an instant, the Jacobian matrix (which
for a potential mapping coincides with the Hessian of the potential) becomes
degenerate; then, generically, the inverse
mapping ceases to exist, i.e., the property of the $(t,\tau)$-mapping to be
bijective gets lost (in the cosmological context this amounts to {\em shell
crossing} leading to {\em multi-streaming}).

Below, we establish criteria for omni-potentiality. In Sec.~\ref{sub:commu} we
prove that a flow is omni-potential, whenever the Hessians of the potentials
$\Phi(\q,t;\tau)$, calculated at any two points of a trajectory, commute.
In Sec.~\ref{sub:app1} we present an equivalent condition for the commutativity
of Hessians: along any trajectory, at any time $t$,
the Hessian ${\cal H}(t)={\cal H}(\Phi(\q; 0,t))$ and its time
derivative should commute; this is used to show that
omni-potentiality is equivalent to potentiality of both the Lagrangian
and the Eulerian flow velocities (plus the convexity constraints). In
Sec.~\ref{sub:linea} we discuss simple examples of omni-potentiality:
Zeldovich and Zeldovich-type flows. Finally,
in Sec.~\ref{sub:eqn2d} we derive a partial differential equation
for the potential of a two-dimensional omni-potential flow. It
states that a suitable expression constructed from second-order derivatives
depends only on the Lagrangian coordinates, but not on time.
In other words, the Hessian of the potential possesses an {\em invariant},
whose value depends only on the trajectory. We construct similar invariants
in higher-dimensional spaces in the Appendix.

\subsection{Commutation of Hessians of the potential}
\label{sub:commu}

The semigroup associativity \eqref{associative} involves the composition
of two potential maps. We shall show that, in general, such a composition is
potential if and only if the Hessians commute. Basically this stems from the
well-known theorem that the product of two symmetric matrices is symmetric only
when they commute. The problem we are now addressing is illustrated
in Fig.~\ref{f:fig1}, which sketches the action of three mappings along the
same trajectory. We observe that the $(t,\tau)$-mapping is the composition of
the inverse of the $(t_0,t)$-mapping with the $(t_0,\tau)$-mapping. We shall
now show that its potentiality is equivalent to the commutation of the Hessians
of the $(t_0,t)$-mapping and of the $(t_0,\tau)$-mapping.

Specifically, we assume that, for any times $t$ and $\tau$, such that $0\le t\le\tau\le T$,
the $(t,\tau)$-mapping \rf{potmap} of $\R^d$ into itself is a bijection
that, together with its inverse, is smooth (i.e., has as many continuous
derivatives, as we might need). We denote
$$\Phi_1(\q)\equiv\Phi(\q,t_0;t),\qquad\Phi_2(\q)\equiv\Phi(\q,t_0;\tau),$$
where $\Phi(\q,t_0;t)$ and $\Phi(\q,t_0;\tau)$ are the potentials
of the $(t_0,t)$-mapping and the $(t_0,\tau)$-mapping, respectively.
The required potentiality of the $(t,\tau)$-mapping implies that the Jacobian
matrix $\left\|{\partial{\bm\xi}\over\partial\x}\right\|$ is symmetric:
$$\partial_{x_j}\xi_i=\partial_{x_i}\xi_j$$
for any pair of indices $1\le i,j\le d$. The converse is also true, at least
locally in space. By the chain rule,
\begin{eqnarray*}
{\cal H}_{mn}(\Phi_2)&=&\partial_{q_n}\xi_m\\
&=&\sum_{k=1}^d\partial_{x_k}\xi_m\,\partial_{q_n}x_k\\
&=&\sum_{k=1}^d\partial_{x_k}\xi_m\,{\cal H}_{kn}(\Phi_1).
\end{eqnarray*}
Therefore,
\begin{equation}
\left\|{\partial{\bm\xi}\over\partial\x}\right\|={\cal H}(\Phi_2)
{\cal H}^{-1}(\Phi_1),
\label{commute}\end{equation}
where ${\cal H}^{-1}(\Phi_1)$ denotes the matrix inverse to
${\cal H}(\Phi_1)$. For the r.h.s.~of \eqref{commute} to be a symmetric
matrix, the matrices ${\cal H}(\Phi_2)$ and ${\cal H}^{-1}(\Phi_1)$ must
commute, which is equivalent to the commutation of the two Hessians
${\cal H}(\Phi_2)$ and ${\cal H}(\Phi_1)$.

Reversing the arguments, we establish that commutation of
Hessians $\Phi(\q,t_0;t)$ along a trajectory,
$${\cal H}(\Phi(\q,t_0;t)){\cal H}(\Phi(\q,t_0;\tau))\!=\!
{\cal H}(\Phi(\q,t_0;\tau)){\cal H}(\Phi(\q,t_0;t)),$$
at any times $t$ and $\tau$ is \textit{necessary} for the mapping
${\bm\xi}(\x)$ to be potential. In particular, commutation of Hessians
of the potential $\Phi(\q,0;t)$ of the Lagrangian map along each trajectory
(for each fixed $\q$), together with invertibility or convexity, is equivalent
to omni-potentiality.

By the theorem on codiagonalizability of commuting (real) symmetric matrices
(see, e.g., Ref.~\cite{hornjohnson}, pp.~50--51) the equivalent condition is
that the Hessians of the potential $\Phi(\q,t_0;t)$ calculated at different
times $t$ for the same coordinate $\q$ and the same $t_0$ are
{\em codiagonalizable}, i.e., can be transformed into the diagonal form using
the same unitary matrix. In other words, along any trajectory, only
the eigenvalues of the Hessian of the potential but not the eigendirections
are allowed to vary.

\begin{figure}
\includegraphics[width=85mm]{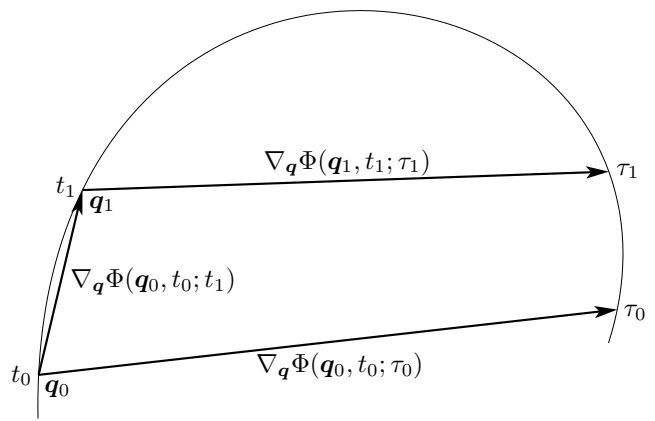}
\caption{Commutation of Hessians of the potential along a trajectory.
A sketch of a trajectory and flow-induced mappings
from times $t_0$ and $t_1$ to times $\tau_0$ and $\tau_1$.}
\label{f:fig2}
\end{figure}

So far, we have shown the commutation --- along a given trajectory ---
of the Hessians of the potentials of the $(t_0,t)$-mapping and
the $(t_0,\tau)$-mapping for the same starting time $t_0$ (e.g., for the
$(t_0,\tau_0)$-mapping and the $(t_0,t_1)$-mapping of Fig.~\ref{f:fig2}).
A similar argument proves the commutation of the Hessians of the
potentials of two mappings, such that the ending time of one of them coincides
with the starting time of the second one (e.g., for $(t_0,t_1)$-mapping and
the $(t_1,\tau_1)$-mapping). Combining these two results and relying again
on the theorem on codiagonalizability of symmetric commuting matrices,
we find that the Hessians of the potentials of the $(t_0,\tau_0)$-mapping and
of the $(t_1,\tau_1)$-mapping commute. Thus we have established that, along any
given trajectory, the Hessians associated with arbitrary pairs of times commute.

\subsection{Commutation of Hessians and their time derivatives:
bi-potential velocities}
\label{sub:app1}

Here we give an alternative formulation of omni-potentiality in terms
of commutation of Hessians and their time derivatives. We need some
preparatory material regarding $d\times d$ symmetric matrices with
smooth time dependence. Suppose at any two times $t$ and $t'$ they commute:
\begin{equation}
H(t)H(t')=H(t')H(t).
\label{tcomm}\end{equation}
Differentiating this equation in $t'$ and letting $t'=t$, we find that
at any time $H(t)$ commutes with its time derivative $\dot{H}(t)$:
\begin{equation}
H(t)\dot{H}(t)=\dot{H}(t)H(t).
\label{dcomm}\end{equation}
We shall show now that the converse is also true: if at any time $t$ (i)
relation \rf{dcomm} is satisfied, and (ii) all eigenvalues $\lambda_i$
of the symmetric matrix $H(t)$ are distinct, then \rf{tcomm} holds true
for any times $t$ and $t_1$.

Since $H(t)$ is symmetric, it can be expressed as
\begin{equation}
H(t)=U^t(t)\Lambda(t)U(t),
\label{dnf}\end{equation}
where $U$ is a unitary $d\times d$ matrix, and $\Lambda$ is diagonal.
Consider the identity $U(t)U^t(t)=I$, where $I$ is
the identity matrix; differentiating it in time yields
$$U\dot{U}^t=-\dot{U}U^t=-\left(U\dot{U}^t\right)^t.$$
Thus $X\equiv U\dot{U}^t$ is an antisymmetric matrix.

We substitute \rf{dnf} into \rf{dcomm}, and multiply it by $U$
on the left and by $U^t$ on the right, obtaining
$$\Lambda X\Lambda-\Lambda^2X=X\Lambda^2-\Lambda X\Lambda,$$
i.e., the matrix $X\Lambda-\Lambda X$ commutes with $\Lambda$.
This matrix is symmetric due to the antisymmetry of $X$.
By the theorem on codiagonalizability of commuting symmetric matrices,
the matrices $X\Lambda-\Lambda X$ and $\Lambda$ are thus
simultaneously diagonalizable; but $\Lambda$ is already diagonal, and hence
so is $X\Lambda-\Lambda X$. The entries of the latter matrix are
$(\lambda_j-\lambda_i)X_{ij}$, and therefore $X_{ij}=0$ for all
$i\ne j$ (we use here the condition that all eigenvalues of $H$
are distinct). The antisymmetry of $X$ implies that all diagonal
entries of $X$ also vanish, and thus $X=0$. Therefore,
$$\dot{U}=-XU$$
vanishes. In other words, variation of $H$ in time consists
solely of variation of its eigenvalues $\lambda_i$. This implies \rf{tcomm}.

The restriction that all the eigenvalues be distinct does not significantly
affect the generality of the statement: if they are not distinct at some
isolated times, the relation $\dot{U}=0$ remains satisfied at these times
by continuity.

Returning to the problem of omni-potentiality, we now take
$$H(t)={\cal H}(\Phi(\q,0;t))\equiv{\cal H}.$$
This Hessian $\cal H$ is the Jacobian of the Lagrangian map $\nabla\Phi(\q,0;t)$.
By the above result, omni-potentiality is equivalent to the commutation,
at any time, of the Jacobian and its time derivative.

We have seen earlier that omni-potential flow has a velocity which
is bi-potential, i.e., potential in both Eulerian and Lagrangian coordinates.
The statement just derived allows to prove the converse, namely that
a flow with a bi-potential velocity (and some convexity requirements)
is omni-potential. Let us denote by ${\bf v}^{\rm L}(\q,t)$ and
${\bf v}^{\rm E}(\x,t)$ the Lagrangian and Eulerian velocity, respectively.
Since $\x(\q,t)$ is the Lagrangian map, we obviously have
\begin{equation}
{\bf v}^{\rm E}(\x,t)={\bf v}^{\rm L}(\q(\x,t),t),
\label{vlageul}\end{equation}
where $\q(\x,t)$ is the inverse Lagrangian map, whose Jacobian is
${\cal H}^{-1}$. We now calculate the Eulerian velocity gradient,
using \eqref{vlageul}. By the chain rule, for any $i$ and $j$, we have:
\begin{eqnarray*}
\partial_{x_i}v^{\rm E}_j(\x,t)&=&\sum_{m=1}^d
({\cal H}^{-1})_{im}\,\partial_{q_m}v^{\rm L}_j(\q,t)\\
&=&\sum_{m=1}^d({\cal H}^{-1})_{im}\,\partial_{q_mq_j}\dot{\Phi}(\q,0;t)\\
&=&\sum_{m=1}^d({\cal H}^{-1})_{im}\dot{\cal H}_{mj}.
\end{eqnarray*}

Thus, the Eulerian velocity gradient is the product of the matrices
${\cal H}^{-1}$ and $\dot{\cal H}$. For the Eulerian velocity to
be potential, it is necessary and (locally) sufficient that this
product be a symmetric matrix. The commutation of the symmetric
matrices ${\cal H}^{-1}$ and $\dot{\cal H}$ is equivalent to the commutation of
${\cal H}$ and $\dot{\cal H}$. Equivalence to omni-potentiality follows from
the statement derived above.

\subsection{Zeldovich and Zeldovich-type flows}
\label{sub:linea}

In the Zeldovich approximation each particle keeps its initial velocity unaltered
in the course of time, and hence particles move along straight lines (at least
before multi-streaming occurs). The Lagrangian map at time $t$ is
$$\q\mapsto\x=\nabla_{\q}\left({|\q|^2\over2}+t\varphi_0(\q)\right),$$
where $\varphi_0(\q)$ is the velocity potential, prescribed at $t=0$.
The Hessian of this map is $I+t{\cal H}(\varphi_0)$, where $I$ is the identity
matrix and the matrix ${\cal H}(\varphi_0)$ is the Hessian of the initial
potential. For a given $\q$, the eigendirections of the associated Hessian are
those of the Hessian of the initial velocity potential. Clearly, all these
Hessians commute and, by the results of Sec.~\ref{sub:commu},
such a flow is omni-potential.

More general examples of omni-potential flows can be constructed by performing
an arbitrary nonlinear transformation of the time and by time-dependent zooming
of space. In space of any dimension $d\ge 2$, consider the flows defined by
the potentials
\begin{equation}
\Phi(\q,0;t)=\mu(t){|\q|^2\over 2}+\eta(t)\varphi_0(\q),
\label{genZel}\end{equation}
where $\mu(t)$ and $\eta(t)$ are arbitrary functions of time. Clearly,
these are again omni-potential.

In general, the trajectories associated with \eqref{genZel} are not straight
lines. However, if we look at them with a time-dependent magnifying glass which
applies a zooming factor $1/\mu(t)$, they become
straight. Furthermore, if we introduce a new time variable
$t'=\eta(t)/\mu(t)$, particles move again with a constant velocity.
Hence, the flows defined by \eqref{genZel} are trivial generalizations
of Zeldovich flows, and will here be called \textit{Zeldovich-type} flows.

Our goal is to find omni-potential flows that are not of this type.

\subsection{A linear second-order PDE for two-dimensional omni-potential flow}
\label{sub:eqn2d}

We derive here a differential equation for the potential of a two-dimensional
omni-potential flow. It turns out to be a linear second-order PDE.

Consider a symmetric $2\times2$ matrix $H$.
Suppose its eigenvector associated with eigenvalue $\lambda$
makes angle $\theta$ with the cartesian axis $q_1$. Thus,
$$H_{11}\cos\theta+H_{12}\sin\theta=\lambda\cos\theta,$$
$$H_{12}\cos\theta+H_{22}\sin\theta=\lambda\sin\theta.$$
In order to eliminate the eigenvalue $\lambda$, we multiply the first of these
equations by $\sin\theta$ and the second one by $\cos\theta$. Subtracting
afterwards the second equation from the first one, we obtain
\begin{equation}
{H_{11}-H_{22}\over H_{12}}=\cot2\theta.
\label{hess2d}\end{equation}
Prescribing the r.h.s.~of \eqref{hess2d} uniquely defines the orthogonal
frame of the two eigendirections. (The values of $\cot2\theta$ define the angle
$\theta$ modulo $\pi/2$; however, changing $\theta\to\theta+\pi/2$ swaps the
eigendirections, but does not affect the set of eigendirections.)

In an omni-potential flow, the eigendirections of the Hessians
of the $(0,t)$-potentials should depend only on the Lagrangian position $\q$
and not on the time $t$. Let $\Phi(\q,t)$ be a two-dimensional omni-potential
flow and let us denote by $g(\q)$ the common value of $\cot2\theta$ along
the particle trajectory emanating from $\q$. It then follows from \eqref{hess2d} that
\begin{equation}
(\partial^2_{q_1q_1}-\partial^2_{q_2q_2})\Phi=g(\q)\,\partial^2_{q_1q_2}\Phi
\label{gen2d}.\end{equation}

The search for two-dimensional omni-potential flow has thus been reduced to
finding solutions to \eqref{gen2d} for suitably prescribed functions $g(\q)$.

\section{Examples of omni-potential flows in two and three dimensions}
\label{s:algebra}

The main question that we address in this paper is whether omni-potential
flows exist that are not of Zeldovich type. In this section we give
a positive reply to this question both in the two- and three-dimensional spaces
by providing explicit examples of polynomial potentials for mappings induced
by such flows.

An example of an omni-potential flow in a space of arbitrary dimension is provided
by spherically-symmetric potentials of the form $\Phi(|\q|,t)$: a simple
calculation reveals the commutation of Hessians of such potentials, calculated
at different times at different points of a trajectory. This example
shares with Zeldovich flows the property that the trajectories are
straight lines --- in this case, in the radial direction. We would like
to construct less symmetric examples.

\subsection{Particular examples of two-dimensional omni-potential flow}
\label{sub:ex2d}

In $\R^2$, the problem of finding omni-potential flows has been recast
into the form of the partial differential equation \rf{gen2d} with the
initial condition $|\q|^2/2$ (which generates the identity map). We can
therefore construct examples of two-dimensional omni-potential flows
by finding different solutions to \rf{gen2d} for a given function
$g(\q)$ in the r.h.s. By linearity, any linear combination of such
solutions with time-dependent coefficients is also a solution to
\rf{gen2d}. For example, if $\Phi_1(\q)$ and $\Phi_2(\q)$ are two
sufficiently smooth independent solutions that are also independent of
$|\q|^2/2$ (which is always a solution to \rf{gen2d}), then the flow
with the potential $|\q|^2/2+\alpha_1(t)\Phi_1(\q)+\alpha_2(t)\Phi_2(\q)$
is omni-potential, and is typically not of Zeldovich type; for this,
the functions of time $\alpha_1$ and $\alpha_2$ must be linearly independent
and sufficiently small, so as not to spoil the convexity stemming
from the $|\q|^2/2$ term.

When $g(\q)$ is a ratio of homogeneous polynomials of $\q$ (say, of the same
degree $m$), solutions to \rf{gen2d} can be obtained by a purely algebraic
method. A solution can be sought
in the form of a homogeneous polynomial, $p^{(2)}_n(\q)$, of degree $n\ge m+2$;
then \rf{gen2d} reduces to a system of $m+n-1$ equations for the coefficients of
$p^{(2)}_n(\q)$ and $g(\q)$. (In what follows, $p^{(d)}_n$ denotes certain
homogeneous polynomials of degree $n$ defined in $\R^d$.) The function $g(\q)$
involves $2m+1$ independent coefficients (since the numerator and denominator
can be multiplied by any constant without changing $g(\q)$), and the polynomial
$p^{(2)}_n(\q)$ involves $n$ independent coefficients (since a solution
to \rf{gen2d} can be multiplied by any constant without yielding a new
independent solution). Comparison of the number of equations, $m+n-1$, with
the total number of the unknown coefficients, $2m+n+1$, suggests that we can
construct a family of such solutions, parameterized by $m+2$ coefficients
of $g(\q)$. However, the system of equations for the coefficients is,
in general, nonlinear, and hence its solvability cannot be established just
by counting the numbers of the unknowns and equations. When $g(\q)$ is the ratio
of linear functions, the equations for the coefficients of $p^{(2)}_n(\q)$ are
linear, and can be solved for any prescribed coefficients of~$g(\q)$.

Since the potential $\Phi$ is required to be a convex function on the entire
plane $\R^2$, we start by seeking homogeneous polynomials $p^{(2)}_n(\q)$
involving only even powers of $q_1$ and $q_2$. An instance of a solvable linear
system of equations yielding the coefficients of such polynomials is obtained for
\begin{equation}
g(\q)={aq_1^2-bq_2^2\over q_1q_2},
\label{gexa}\end{equation}
where the coefficients $a$ and $b$ may take arbitrary preset values. For such
$g(\q)$, a homogeneous polynomial of degree $2k$, $k\ge2$, solving \rf{gen2d}, is
\begin{eqnarray}
&&\hspace*{-8mm}p^{(2)}_{2k}(q_1,q_2)=\sum_{i=0}^k\left(\,\prod_{j=0}^{i-1}(2k-1+2j(a-1))\right.\nonumber\\
&&\hspace*{-8mm}\times\left.\prod_{j=0}^{k-1-i}(2k-1+2j(b-1))\right){k!\ q_1^{2i}q_2^{2(k-i)}\over i!(k-i)!(2k-1)}.
\label{poly2k}\end{eqnarray}
In particular, the first low-degree polynomial solutions are:
\begin{eqnarray}
p^{(2)}_4(q_1,q_2)&=&(2a+1)q_1^4+6q_1^2q_2^2+(2b+1)q_2^4,\label{p4}\\
p^{(2)}_6(q_1,q_2)&=&(4a+1)(2a+3)q_1^6+15(2a+3)q_1^4q_2^2\nonumber\\
&&\hspace*{-1cm}+15(2b+3)q_1^2q_2^4+(4b+1)(2b+3)q_2^6.\label{p6}
\end{eqnarray}
As one can see, the polynomial \rf{poly2k} vanishes identically for integer
$j\ge1$ and $\widehat{j}\ge1$ such that $j+\widehat{j}\le k-1$, and
\begin{equation}
\widehat{a}=1-{2k-1\over2\widehat{j}},\qquad\widehat{b}=1-{2k-1\over2j}.
\label{param}\end{equation}
For these isolated values in the plane of parameters, there exist nevertheless
two independent solutions, namely
$${\partial\over\partial a}\left.p^{(2)}_{2k}\right|_{a=\widehat{a},\,b=\widehat{b}}
\qquad\hbox{and}\qquad
{\partial\over\partial b}\left.p^{(2)}_{2k}\right|_{a=\widehat{a},\,b=\widehat{b}}.$$
This can be easily seen by differentiating \rf{gen2d} in $a$ and $b$ and
substituting the parameter values \rf{param}.

Clearly, $p^{(2)}_{2k}(\q)$ is convex, if all coefficients in \rf{poly2k} are
positive, i.e., if
\begin{equation}
\min(a,b)\ge-1/(2k-2).
\label{minabe}\end{equation}
Thus, the potentials
\begin{eqnarray}
\Phi(\q,t)=\mu_2(t){|\q|^2\over 2}+\sum_{k\ge2}\mu_{2k}(t)p^{(2)}_{2k}(q_1,q_2)
\label{exa2d}\end{eqnarray}
are convex for $\min(a,b)\ge0$, if in addition all $\mu_{2k}(t)$ are non-negative
(this condition is sufficient, but not necessary) and tend to zero fast enough
to guarantee convergence of the series at any point $\q$ and termwise differentiability
of \rf{exa2d} in the spatial variables. If the sum \rf{exa2d} is finite
and the maximum degree of the polynomials involved is $2K$, then the potential
is convex for $\min(a,b)\ge-1/(2K-2)$. The initial condition is satisfied
provided $\mu_2(0)=1$ and $\mu_{2k}(0)=0$ for all $k>1$. The convex potentials
\rf{exa2d} satisfy all requirements for omni-potential flows in the plane, and
are not of Zeldovich type.

So far, we have considered only even-degree homogeneous polynomial solutions or
linear combinations thereof. Is an admixture of odd-degree homogeneous
polynomials permitted? If such an odd-degree addition has a degree higher than
that of the highest even-degree homogenous polynomial comprised in the
solution, then convexity in the whole plane is ruled out. However, in a finite
domain, convexity need not be lost, provided the odd polynomial is scaled by
a sufficiently small factor. This is precisely what happens when $g(\q)$ is given by
\rf{gexa}: a homogeneous polynomial of odd degree $2k+1$ can be a solution to
\rf{gen2d} only for $k\ge2$ and
\begin{equation}
a=b=-1/(k-1).
\label{minabo}\end{equation}
The solution for these parameter values is
$$p^{(2)}_{2k+1}(q_1,q_2)=c_1\pf_k(q_1,q_2)+c_2\pf_k(q_2,q_1),$$
where
\begin{eqnarray*}
\pf_k(q_1,q_2)\!=\!\sum_{i=0}^{k-1}&&\hspace*{-4mm}\left(\,
\prod_{j=0}^{i-1}{(2(k-j)+1)(k-1-j)\over(j+1)(2j-1)}\right)\\
&&\hspace*{-3mm}\times\,q_1^{2i}q_2^{2(k-i)+1},
\end{eqnarray*}
and $c_1$ and $c_2$ are arbitrary constants.
For instance, for $a=b=-1$ fifth-degree polynomial solutions are
$$p^{(2)}_5(q_1,q_2)=c_1q_1^5-5c_1q_1^3q_2^2-5c_2q_1^2q_2^3+c_2q_2^5.$$
Comparison of conditions \rf{minabe} and \rf{minabo} shows that
the degree of an odd-degree polynomial solution for $g(\q)$ defined by \rf{gexa}
is higher than the degree of any convex even-degree polynomials existing
for the chosen $g(\q)$.

\subsection{Examples of omni-potential flow in dimension $d\ge 3$}
\label{sub:ex3d}

The approach that has been applied in the previous subsection for construction
of an example in dimension two cannot be immediately generalized
to higher-dimensional spaces: while in $\R^2$ a single invariant determines
whether two symmetric matrices are codiagonalizable and hence equation \rf{hess2d}
fixes the set of eigendirections of a symmetric matrix $H$, in $\R^3$ at least
three such invariants must be considered simultaneously (see the Appendix).
In dimension three, equations \rf{inv3d1}--\rf{inv3d2} applied for the entries
of the Hessian of an unknown potential give rise to three partial differential
equations in the potential,
\begin{eqnarray}
&&{\partial^2_{q_1,q_3}\Phi\over\partial^2_{q_1,q_2}\Phi}\!\left(
{(\partial^2_{q_1,q_1}-\partial^2_{q_2,q_2})\Phi\over\partial^2_{q_1,q_2}\Phi}
+{\partial^2_{q_2,q_3}\Phi\over\partial^2_{q_1,q_3}\Phi}
-{\partial^2_{q_1,q_3}\Phi\over\partial^2_{q_2,q_3}\Phi}\right)\nonumber\\
&=&\!\!\!\left(g_1(\q)+{(\partial^2_{q_1,q_1}-\partial^2_{q_2,q_2})\Phi\over\partial^2_{q_1,q_2}\Phi}\right)\nonumber\\
&\times&\!\!\!\left(
{(\partial^2_{q_2,q_2}-\partial^2_{q_3,q_3})\Phi\over\partial^2_{q_2,q_3}\Phi}
+{\partial^2_{q_1,q_3}\Phi\over\partial^2_{q_1,q_2}\Phi}
-{\partial^2_{q_1,q_2}\Phi\over\partial^2_{q_1,q_3}\Phi}\right),\label{pde3d1}\\
&&{\partial^2_{q_2,q_3}\Phi\over\partial^2_{q_1,q_3}\Phi}
\left(g_1(\q)+{(\partial^2_{q_1,q_1}-\partial^2_{q_2,q_2})\Phi\over\partial^2_{q_1,q_2}\Phi}\right)\nonumber\\
&=&\!\!g_2(\q)-g_3(\q){(\partial^2_{q_1,q_1}-\partial^2_{q_2,q_2})\Phi\over\partial^2_{q_1,q_2}\Phi},\label{pde3d2}\\
&&g_3(\q)\left(
{(\partial^2_{q_2,q_2}-\partial^2_{q_3,q_3})\Phi\over\partial^2_{q_2,q_3}\Phi}
+{\partial^2_{q_1,q_3}\Phi\over\partial^2_{q_1,q_2}\Phi}
-{\partial^2_{q_1,q_2}\Phi\over\partial^2_{q_1,q_3}\Phi}\right)\nonumber\\
&=&\!\!{(\partial^2_{q_3,q_3}-\partial^2_{q_1,q_1})\Phi\over\partial^2_{q_1,q_3}\Phi}
+{\partial^2_{q_1,q_2}\Phi\over\partial^2_{q_2,q_3}\Phi}
-{\partial^2_{q_2,q_3}\Phi\over\partial^2_{q_1,q_2}\Phi},
\label{pde3d3}\end{eqnarray}
that must be satisfied simultaneously. Here the time-independent
quantities $g_k(\q)$ are related to the invariants $\gamma^{(3,k)}_{21}$
introduced in the Appendix:
\begin{eqnarray*}
g_1(\q)&=&\gamma^{(3,1)}_{21}+\gamma^{(3,3)}_{21},\\
g_2(\q)&=&\gamma^{(3,2)}_{21}+1,\\
g_3(\q)&=&\gamma^{(3,3)}_{21}.
\end{eqnarray*}
The invariants are combinations of the entries of a symmetric matrix,
which depend only on the set of eigendirections and not on the eigenvalues.
The algebraic theory of the invariants presented in the Appendix does not
take into account the properties of the Hessian, stemming from the specific
structure of its entries (for instance, each row and a column of the Hessian is
a gradient, which implies certain differential relations between the entries).
It is unclear how to prescribe $g_k(\q)$, taking into account these additional
constraints, for the equations to have at least two distinct solutions.

Because of this difficulty, instead of considering the invariants and solving
equations \rf{pde3d1}--\rf{pde3d3}, we exploit the fact that omni-potentiality
of flows amounts to commutation of the Hessians of the various
$(t_1,t_2)$-mappings along any trajectory (see Sec.~\ref{sub:commu}). We shall
construct our examples in $\R^d$ using the following strategy. The potential
is sought in the form of a linear combination of ``building blocks'' with
time-dependent coefficients. One of the building blocks is prescribed; we
take it to be a homogeneous polynomial, $p^{(d)}_m(\q)$, of degree $m$. All
the other building blocks are then required to have their Hessians commuting
with that of the prescribed building block. The function $|\q|^2$, whose
Hessian is the identity matrix, constitutes a trivial solution. We can try
finding other building blocks in the form of homogeneous polynomials
$p^{(d)}_n(\q)$ of some degree $n$. If we succeed, it is easy to show that
any linear combination of such building blocks (with the convexity restriction)
will define an omni-potential flow. We seek such polynomials by requiring the
vanishing of the non-diagonal entries of the commutator of the two Hessians, viz.
\begin{equation}
C(p^{(d)}_m,p^{(d)}_n)\equiv{\cal H}(p^{(d)}_m){\cal H}(p^{(d)}_n)-{\cal H}(p^{(d)}_n){\cal H}(p^{(d)}_m).
\label{commr}\end{equation}

Unfortunately, in general this strategy does not work, as now explained.
A homogeneous polynomial of degree $n$ has
$$(n+d-1)!\over n!(d-1)!$$
coefficients. Non-diagonal entries of the commutator $C$ are homogeneous
polynomials of degree $m+n-4$. The commutator is antisymmetric (recall that
the Hessians are symmetric matrices), hence we have to consider the $d(d-1)/2$
non-diagonal entries of $C$. Thus, in general, we have to solve
$$d(m+n+d-5)!\over2(m+n-4)!(d-2)!$$
equations, a number which is easily seen to exceed the number of coefficients,
$${(m+d-1)!\over m!(d-1)!}+{(n+d-1)!\over n!(d-1)!}.$$
So, the problem is overdetermined.

Nevertheless, potentials having all the required properties can be constructed
in $\R^d$ ($d\ge2$), if all the building blocks are restricted to be
homogeneous polynomials symmetric in all their arguments (i.e., invariant
under any permutation of the spatial variables $q_i\leftrightarrow q_j$). Such
building blocks have the following significant advantage: it suffices to
consider any of the polynomial equations arising from non-diagonal entries
of the commutator \rf{commr} (referred to as ``commutator equations'') -- all
these equations are equivalent by virtue of the symmetry.
In what follows, we implement this ``symmetric building block strategy''
in two cases, in $\R^d$ for $d\ge 3$ with just one unknown building block,
and in $\R^3$ with infinitely many ones.

Now, we focus on the symmetric polynomials
\begin{eqnarray}
p^{(d)}_4(\q)&=&\sum_{i=1}^dq_i^4+\widetilde{c}\sum_{i=2}^d\sum_{j=1}^{i-1}q_i^2q_j^2,
\label{p4gen}\\
p^{(d)}_6(\q)&=&\sum_{i=1}^dq_i^6+\widetilde{a}\sum_{i=1}^d\sum_{j=1}^dq_i^4q_j^2\nonumber\\
&&+\,\widetilde{b}\sum_{1\le i<j<k\le d}q_i^2q_j^2q_k^2.
\label{p6gen}\end{eqnarray}
(When $d=3$, the last sum in $p^{(d)}_6(\q)$ reduces to a single term
$\widetilde{b}q_1^2q_2^2q_3^2$.) We consider polynomials involving only even
powers of the spatial variables $q_j$, because we seek solutions that are
convex functions.

To be specific, we consider the commutator equation
$C_{12}(p^{(d)}_6,p^{(d)}_4)=0$ in $\R^d$, $d\ge3$. The l.h.s.~is a polynomial
of degree 6. Since in $p^{(d)}_4$ and $p^{(d)}_6$ any power of $q_1$ and $q_2$
is even, $C_{12}$ is proportional to $q_1q_2$, and every variable enters
into the polynomial $C_{12}/(q_1q_2)$ only in an even power.
Since both potentials are symmetric in $q_1$ and $q_2$, $C_{12}=0$
for $q_1=q_2$, and hence $C_{12}$ is divisible by
$q_1^2-q_2^2$. The polynomial $C_{12}/(q_1q_2(q_1^2-q_2^2))$
is of the second degree; it is thus just a sum of $q_j^2$ with certain
coefficients. Because the potentials are symmetric in $q_j$, it has the form
$$\alpha_1(q_1^2+q_2^2)+\alpha_2\sum_{j=3}^dq_j^2.$$
Hence, we have three independent parameters, $\widetilde{a}$, $\widetilde{b}$ and
$\widetilde{c}$, and two equations to satisfy. Calculating the coefficients
$\alpha_1$ and $\alpha_2$ and letting them vanish, we find
that the Hessians of $p^{(d)}_4(\q)$ and $p^{(d)}_6(\q)$ commute for
\begin{eqnarray}
\widetilde{a}&=&15\,\widetilde{c}/(12-\widetilde{c}),\label{p6a}\\
\widetilde{b}&=&75\,\widetilde{c}^{\,2}/((12-\widetilde{c})(3+\widetilde{c})).
\label{p6b}\end{eqnarray}
Thus, the potential
\begin{equation}
\Phi(\q,t)=\mu_2(t){|\q|^2\over 2}+\mu_4(t)p^{(d)}_4(\q)+\mu_6(t)p^{(d)}_6(\q)
\label{polydd}\end{equation}
defines a non-Zeldovich-type omni-potential flow in $\R^d$ ($d\ge3$).
Polynomials $p^{(d)}_4(\q)$ and $p^{(d)}_6(\q)$ are convex provided
$0\le\widetilde{c}<12$; hence, if all $\mu_i(t)\ge 0$, potential
\rf{polydd} is convex for $\widetilde{c}$ from this interval.
For $\widetilde{c}\ne 2$, $p^{(d)}_4(\q)$ and $p^{(d)}_6$ do not possess
spherical symmetry, and hence the potential \rf{polydd} is not spherically
symmetric. Restrictions of $p^{(d)}_4(\q)$ and $p^{(d)}_6(\q)$ onto the plane
$q_3=...=q_d=0$ coincide with the polynomials $p^{(2)}_4(q_1,q_2)$ \rf{p4} and
$p^{(2)}_6(q_1,q_2)$ \rf{p6} for $a=b=(6-\widetilde{c})/(2\widetilde{c})$.

Henceforth, for the sake of simplicity, we assume that the problem is
three-dimensional. In the remainder of the subsection we shall implement
the building block strategy with unknown blocks that are arbitrary even-degree
polynomials $p^{(3)}_{2n}(\q)$ for $n>2$. The polynomial $p^{(d)}_4(\q)$
\rf{p4gen} for $d=3$ remains our prescribed building block. By the theorem
on codiagonalizability of symmetric matrices, commutation with the Hessian
of $p^{(3)}_4(\q)$ implies, that the Hessians of any two polynomials from this
family commute. (This is taking place generically, i.e., at those points
in $\R^3$, where the Hessian of $p^{(3)}_4(\q)$ does not possess equal
eigenvalues; at non-generic points the commutation follows from continuity
of the Hessians and the commutation at generic points, which are present
in any neighborhood of a non-generic point.)

The polynomial
\begin{equation}
p^{(3)}_{2n}(\q)=\sum_{\mbox{$\scriptsize i,j,k\ge0\atop i+j+k=n$}}
\widetilde{a}_{i,j,k}\,q_1^{2i}q_2^{2j}q_3^{2k}
\label{p2n}\end{equation}
is symmetric, whenever
\begin{equation}
\mbox{$\widetilde{a}_{i,j,k}$ does not depend on the order of subscripts $i,j,k$.}
\label{p2ncoe}\end{equation}
Straightforward algebra yields
\begin{eqnarray*}
&&\hspace*{-6mm}C_{12}(p^{(3)}_{2n},p^{(3)}_4)=8q_1q_2\sum_{\mbox{$\scriptsize i,j,k\ge0\atop i+j+k=n$}}
\widetilde{a}_{i,j,k}\,q_1^{2i-2}q_2^{2j-2}q_3^{2k}\\
&&\hspace*{4mm}\times\left(ij(\widetilde{c}-6)(q_1^2-q_2^2)\right.\\
&&\hspace*{4mm}\left.+\,\widetilde{c}\,(-j(2j-1+2k)q_1^2+i(2i-1+2k)q_2^2)\right).
\end{eqnarray*}
Collecting similar terms in this sum, we find that it vanishes as long as
\begin{equation}
\widetilde{a}_{i,j,k}=\widetilde{a}_{i+1,j-1,k}\,\chi_j/\chi_{i+1}
\label{recur}\end{equation}
for any $i,\,j$ and $k$ such that $i+j+k=n$, where we have denoted
$$\chi_m=(\widetilde{c}\,(2n+2-3m)+6(m-1))/m.$$
Relation \rf{recur} can be regarded as a recurrence for coefficients
$\widetilde{a}_{i,j,k}$ for a fixed $k$. For $k=0$, we start the recurrence
assuming $\widetilde{a}_{n,0,0}=1$. This yields
$$\widetilde{a}_{i,n-1,0}={\prod_{m=1}^{n-i}\chi_m\prod_{m=1}^i\chi_m\over\prod_{m=1}^n\chi_m}.$$
We obtain now the starting values for the recurrence \rf{recur} for $k>0$
setting, in view of \rf{p2ncoe},
$$\widetilde{a}_{n-k,0,k}=\widetilde{a}_{n-k,k,0},$$
and find
\begin{equation}
\widetilde{a}_{i,j,k}={\prod_{m=1}^i\chi_m\prod_{m=1}^j\chi_m\prod_{m=1}^k\chi_m\over\prod_{m=1}^n\chi_m}.
\label{coeff2n}\end{equation}
Evidently, coefficients \rf{coeff2n} satisfy the symmetry condition \rf{p2ncoe}.
Hence, \rf{p2n} with the coefficients \rf{coeff2n} is a symmetric homogeneous
polynomial of degree $2n$, whose Hessian commutes with the Hessian of $p^{(3)}_4(\q)$.
The potential
\begin{equation}
\widetilde{\Phi}(\q,t)=\mu_2(t){|\q|^2\over 2}+\sum_{n\ge 2}\mu_{2n}(t)p^{(3)}_{2n}(\q)
\label{xpoly}\end{equation}
defines an omni-potential flow of a non-Zeldovich type in $\R^3$ (provided
the coefficients $\mu_{2n}(t)$ tend to zero sufficiently fast to guarantee convergence
of the series \rf{xpoly} and to allow its termwise differentiation). Since
the polynomial $p^{(3)}_{2n}(\q)$ is convex for
$$0\le\widetilde{c}<6(n-1)/(n-2),$$
the potential \rf{xpoly} is convex if all $\mu_{2n}(t)$ are non-negative and
$0\le\widetilde{c}\le 6$ (or for $0\le\widetilde{c}<6(N-1)/(N-2)$, if all
$\mu_{2n}(t)$ vanish for $n>N$).

Although we have constructed our example without prescribing the invariants,
it might be of interest to calculate them for the solutions that have been
obtained. Straightforward calculations yield the values of the invariants
for the potential $\widetilde\Phi(\q,t)$, for instance,
$$\gamma^{(3,1)}_{21}={(6+3\widetilde{c})q_2^2-(6+\widetilde{c})q_1^2\over2\widetilde{c}q_1q_2}
+{2q_2(q_1^2-q_2^2)\over q_1(q_2^2-q_3^2)},$$
which is consistent with a non-trivial dependence of the eigendirections
on the trajectories (labeled by the Lagrangian coordinates). We note that,
although the solution is symmetric in the spatial coordinates, the symmetry
is lost in the invariant. This stems from the invariant under
consideration being a nonlinear function of the projections of the
eigendirections on the plane $(q_1,q_2)$, and also from the components of the
eigendirections not being invariant under all permutations of coordinates (an eigenvector is
invariant under a permutation of the spatial variables $q_i\leftrightarrow q_j$
provided its $i$th and $j$th components are swapped).

We have used two approaches for constructing examples of omni-potential flow.
In three or more dimensions, we used the building-block strategy in which
the field of eigendirections of the commuting Hessians is characterized by
prescribing one of the blocks (in our construction the polynomial \rf{p4gen}).
The problem of commutation of Hessians then reduces to three linear equations
in the unknown block. These equations must be satisfied simultaneously, and we
found that in general no solution exists for an arbitrary prescribed block.
In two dimensions we have followed another approach, whereby the field of
eigendirections of the commuting Hessians is characterized by prescribing the
set of invariants, from which the field of eigendirections of the commuting
Hessians can be uniquely determined. In $\R^2$ just one such invariant,
$g(\q)$ (see \rf{gen2d}), should be considered. In $\R^3$ one must consider
three invariants, for instance, \rf{inv3d1}--\rf{inv3d2} of the Appendix, giving rise to three
nonlinear equations \rf{pde3d1}--\rf{pde3d3}. As in the former approach, these
equations must be satisfied simultaneously, and hence a solution does not
exist for an arbitrary set of prescribed invariants. Thus, whichever approach
is used for construction of omni-potential flows in $\R^3$, the prescribed
data must be tuned for a solution to exist. In the former approach the
equations are linear and thus simpler, but one is left with just one function
which can be tuned to achieve consistency of the three equations under
consideration; in the latter approach the equations are nonlinear and thus
more involved, but one has the freedom of tuning three a priori independent
scalar functions to gain consistency of the equations. Of course, in three
dimensions, the three invariants of the Hessian of the potential cannot be
prescribed as arbitrary functions. In other words, some conditions
on the invariants must hold for the three equations \rf{pde3d1}--\rf{pde3d3}
to be compatible. Can such conditions on the invariants be expressed in more
explicit form remains an open mathematical problem.

\section{A WKB approach to two-dimensional omni-potentiality}
\label{s:wkb}

So far we have obtained special cases of non-Zeldovich-type omni-potential
flows. How general are they? Can we, for example, in the two-dimensional case
prescribe an arbitrary smooth initial velocity potential $\varphi_0(\q)$ or,
more precisely, the invariant of its Hessian ${\cal H}(\varphi_0(\q))$:
\begin{equation}
g(\q)\equiv\frac{(\duu-\ddd)\varphi_0(\q)}{\dud\varphi_0(\q)}
=\frac{{\cal H}_{11}-{\cal H}_{22}}{{\cal H}_{12}}
\label{gintermsofphi0}\end{equation}
that appears in the general equation \eqref{gen2d}? In this section,
where we use only Lagrangian coordinates, $\du$ and $\dd$ are short
for $\partial_{q_1}$ and $\partial_{q_1}$; similarly, $\duu$,
$\dud$ and $\ddd$ denote the second Lagrangian derivatives.

We shall now show how the construction of non-Zeldovich-type omni-potential flow
with arbitrary invariant function $g(\q)$ can be done, using an idea of Arnold
for solving the linear equation which controls the stability of solution to the
Euler equation \cite{arnold1972}. Rather than trying to find the most general
solution to \eqref{gen2d}, we construct a special short-wavelength
solution through the WKB ansatz
\begin{equation}
\Phi(\q)=\ue^{\ui\kappa S(\q)}\left[A_0(\q)+\frac{1}{\kappa}A_1(\q)
+\frac{1}{\kappa^2}A_2(\q)+\ldots\right]+{\rm cc},
\label{wkb}\end{equation}
where the wavenumber $\kappa$ is taken very large and where cc stands for complex
conjugate (needed because we want real solutions). In WKB parlance, $S(\q)$ is
called the eikonal function and the functions $A_0(\q)$,
$A_1(\q)$, \ldots~are the amplitudes.

To the leading order, O$\left(\kappa^2\right)$, the WKB ansatz turns the linear
second-order PDE \eqref{gen2d} into the following nonlinear first-order PDE:
\begin{equation}
\frac{\left(\du S\right)^2-\left(\dd S\right)^2}{\left(\du S\right)
\left(\dd S\right)}=g(\q).
\label{nonlinearpde}\end{equation}
It is easily checked that \eqref{nonlinearpde} is equivalent to the
statement that, in the leading order, $\nabla_{\q}S(\q)$ is an eigenvector
of the Hessian ${\cal H}(\varphi_0(\q))$. Actually, this can be seen
directly, by an argument which also applies in space dimensions $d$ higher than
two. Assume that the leading WKB term for the potential has a fast spatial
dependence involving the phase factor $\ue^{\ui\kappa S(\q)}$, then the Hessian
will involve in the leading order a matrix factor $-\kappa^2(\partial_i S)(\partial_j S)$.
This is a degenerate matrix with one eigenvector of non-vanishing eigenvalue
in the direction of $\nabla_{\q} S(\q)$; all perpendicular vectors are associated
with the eigenvalue zero, which has multiplicity $d-1$. Omni-potentiality
requires that this degenerate matrix commute with the Hessian of the initial
potential or, equivalently, that $\nabla_{\q}S(\q)$ be an eigenvector
of ${\cal H}(\varphi_0(\q))$.

Returning to the two-dimensional case, we now construct the eikonal
function $S(\q)$. This construction will be done only locally in a
neighborhood $\Omega$, in which the Hessian ${\cal H}(\varphi_0(\q))$ is
sufficiently smooth and its eigenvalues are everywhere distinct. (We
recall that a double eigenvalue is an event of codimension two, which
typically takes place at isolated locations.) Let ${\bm n}^{(1)}(\q)$
and ${\bm n}^{(2)}(\q)$ be two unit eigenvectors of ${\cal H}(\varphi_0(\q))$,
chosen to depend smoothly on $\q$ in $\Omega$. The condition that the gradient
of the eikonal function be parallel to an eigendirection can now be expressed as
\begin{equation}
{\bm n}^{(1)}(\q)\cdot\nabla_{\q}S(\q)=0\quad{\rm or}\quad{\bm n}^{(2)}(\q)\cdot\nabla_{\q}S(\q)=0.
\label{twoeikonal}
\end{equation}
In words, these equations state that either ${\bm n}^{(1)}(\q)$
or ${\bm n}^{(2)}(\q)$ is normal to the level lines of the eikonal function.
Equivalently, the level lines of $S$ are the integral curves defined
by either ${\bm n}^{(2)}(\q)$ or ${\bm n}^{(1)}(\q)$. These form a
set of orthogonal curves. We thus have two classes of solutions.
We can prescribe $S$ arbitrarily on one of these curves ${\cal C}$ and
extend it locally by demanding that it remains constant on all the
curves orthogonal to ${\cal C}$. Note that these orthogonal
curves play here the role of rays in geometrical optics and are thus
conveniently called ``rays''.

Next we write the equations for subleading corrections obtained by substituting
\eqref{wkb} in \eqref{gen2d} and identifying the coefficients of the various
positive and negative powers of the large parameter $\kappa$. We shall
only write the equations appearing at orders $\kappa^1$ and $\kappa^0$
(the higher-order equations have a similar structure). For what follows it is
convenient to use the compact notation introduced by Monge in his theory
of surfaces: $\ph$, $\qh$, $\rh$, $\sh$ and $\th$ stand respectively for
$\du S$, $\dd S$, $\duu S$, $\dud S$ and $\ddd S$. (We added hats to avoid
possible confusions.) Furthermore, we write $g$ for $g(\q)$. The leading-order
equation, \eqref{nonlinearpde}, repeated for convenience, and the two first
subleading equations are:
\begin{eqnarray}
&&\!\!\!\!\!\ph^2 -\qh^2 -g\ph\qh=0, \label{orderm2}\\
&&\!\!\!\!\!(\rh -g\sh -\th)A_0 +2(\ph\du A_0 -\qh\dd A_0)
- g(\ph\dd A_0+\qh\du A_0) \nonumber\\
&&\!\!\!\!\!=0, \label{order1}\\
&&\!\!\!\!\!(\rh-g\sh-\th)A_1 +2(\ph\du A_1-\qh\dd A_1)
- g(\ph\dd A_1+\qh\du A_1) \nonumber\\
&&\!\!\!\!\!-\ui(\duu-\ddd-g\dud)A_0=0.
\label{order0}
\end{eqnarray}

Using \eqref{orderm2} to eliminate the function $g$, we can
rewrite \eqref{order1} and \eqref{order0} as
\begin{eqnarray}
&&\hspace*{-1cm}\qh\du A_0-\ph\dd A_0 +\frac{\ph\qh(\rh-\th)-\sh(\ph^2-\qh^2)}{\ph^2+\qh^2}A_0=0,
\label{order1simpler}\\
&&\hspace*{-1cm}\qh\du A_1-\ph\dd A_1 +\frac{\ph\qh(\rh-\th)-\sh(\ph^2-\qh^2)}{\ph^2+\qh^2}A_1
\nonumber\\
&&\hspace*{-1cm}-\,\ui\frac{\ph\qh}{\ph^2+\qh^2}(\duu-\ddd-g\dud)A_0=0.
\label{order0simpler}
\end{eqnarray}
Equation~\eqref{order1simpler} is a first-order linear homogeneous transport
equation for the amplitude $A_0$ along the rays. It can be integrated starting
from arbitrary non-zero data on any curve orthogonal to the rays.
Equation~\eqref{order0simpler} for the amplitude $A_1$ is of the same sort, except that it has
an inhomogeneous term involving $A_0$. We may thus take vanishing data for $A_1$ on
an arbitrary curve orthogonal to the rays. Higher-order amplitudes satisfy
similar inhomogeneous transport equations.

Now we construct locally in space and time an omni-potential flow having a given
invariant function $g(\q)$. We take the initial potential $\varphi_0(\q)$
arbitrary, but sufficiently smooth. Hence, by the WKB method described
above we can construct a smooth eikonal function $S(\q)$ and smooth amplitude
functions $A_0(\q)$, $A_1(\q)$, \ldots. This, in principle, yields a smooth
solution, $\Phi(\q)$, to \eqref{gen2d}. (We shall not address here the issue
of the convergence of the WKB series \eqref{wkb}.) Because of the imaginary
exponential dependence on $\kappa$, the potential $\Phi(\q)$ has very large
second spatial derivatives O$\left(\kappa^2\right)$ and
has no reason to be convex. However, the following
time-dependent potential defines an omni-potential Lagrangian map:
\begin{equation}
\Phi(\q,t)=\frac{|\q|^2}{2}+t\varphi_0(\q)+f(t)\frac{\epsilon}{\kappa^2}\Phi(\q).
\label{arbitrary2d-omnipot}
\end{equation}
Here $f(t)$ is an arbitrary smooth function of time that vanishes, together
with $\dot f$, at $t=0$. For example, we can take $f(t)=t^2$. For
sufficiently small $t$ and sufficiently small $\epsilon$, the last two terms in
the r.h.s.~of \eqref{arbitrary2d-omnipot} will not spoil the convexity of the
first term. We have thus constructed (locally) omni-potential flows for a
quite arbitrary initial velocity. For large $\kappa$, the trajectories
resulting from \eqref{arbitrary2d-omnipot} differ only minutely from the
straight Zeldovich trajectories associated with the two first terms. However,
these flows are not of Zeldovich type because of the third term in the r.h.s.

If we try to extend the above WKB procedure from two to three
dimensions, we encounter an obstacle already in constructing
the eikonal function $S$. As we have seen, its gradient with respect
to the Lagrangian position ${\q}$ should be everywhere parallel to
an eigendirection of the Hessian ${\cal H}(\varphi_0)$.
Denoting by ${\bm n}(\q)$ a unit eigenvector, taken with
a smooth ${\q}$-dependence, we should then have
\begin{equation}
\nabla_{\q}S(\q)=\mu(\q){\bm n}(\q),
\label{intfactor}
\end{equation}
where $\mu(\q)$ is a scalar function. In other words, the 1-form
${\bm n}(\q)\cdot d\q$ should have an integrating factor (a factor which
makes it an exact 1-form). This is in general possible (locally) in two,
but not in three dimensions.

\section{Cosmological implications}
\label{s:cosmo}

So far our point of view has been kinematical: we constructed omni-potential
flows without any underlying dynamical equations. In cosmology the dynamical
setting is rather well known and discussed for example in
Ref.~\cite{bernardeau}. Let us just recall a few salient points. The most
widely accepted explanation of the large-scale structure seen in galaxy
surveys is that it results from small primordial fluctuations that grew under
gravitational self-interaction of collisionless cold dark matter
particles in an expanding universe.
The evolution of the mostly collisionless matter in the Universe is
described by the Vlasov--Poisson system in the position-velocity phase
space. At early times, i.e., close to the epoch of matter-radiation
decoupling, the expansion of the Universe selects a
single velocity solution at each position rather than a distribution
of velocities. This feature persists until particle crossing
(``shell-crossing'' in the cosmological language), where multi-stream
solutions are developed. As long as multi-streaming is ruled out or is
confined to scales sufficiently small to be neglected, the Vlasov--Poisson
system may be replaced by the Euler-Poisson system. Following the notation
of Ref.~\cite{mnras} and using Eulerian comoving coordinates, $\x$, together
with a time variable $\tau$ based on the amplitude factor of the growing mode
of linear theory, we can write the Euler--Poisson system as
\begin{eqnarray}
\label{comoveuler}
\dtau\v+(\v\cdot\gradx)\v&=&-\frac{3}{2\tau}(\v+\gradx\phig),\\
\label{comovcontinuity}
\dtau\rho+\gradx\cdot(\rho\v) &=& 0,\\
\label{comovpoisson}
\lapx\phig &=& \frac{\rho-1}{\tau}.
\end{eqnarray}
Here, $\v$ is the peculiar velocity, $\rho$ the density (suitably normalized)
and $\phig$ the gravitational potential. As
$\tau\to 0$, to avoid singularities, the density must approach unity
everywhere; thus the distribution of matter is in the leading order uniform as
$\tau\to 0$. Similarly, $\v\to-\gradx\phig$ as $\tau\to 0$;
thus the initial velocity is potential, but
otherwise arbitrary. It then follows from \eqref{comoveuler} that the velocity
stays potential at any later time.

Reconstruction handles the Euler--Poisson system as a two-point
boundary-value problem in which the initial density is prescribed as
uniform and the final (current) density is given by astronomical
observations. This is a mass transport problem whose cost function is
the action associated with the Euler--Poisson equations (see
Refs.~\cite{mnras} and~\cite{loeper}). Unfortunately, because this
action is a rather complicated functional, we do not yet possess
efficient numerical algorithms allowing us to solve this mass transport problem.

The situation simplifies with the Zeldovich approximation, which
amounts to setting to zero the r.h.s.~of \eqref{comoveuler}. The
remaining equations are then (i) the three-dimensional inviscid
Burgers equation, which implies that particles are moving with
constant velocity (in the coordinates here chosen), and (ii) the
continuity equation expressing mass conservation. The action to be
minimized for reconstruction reduces then to its kinetic energy term.
As a consequence, the cost function is just
the mass-weighted integral of the squared displacement of fluid
particles from their initial (Lagrangian) positions $\q$ to their
current (Eulerian) positions $\x$, as required by a theorem of Brenier
\cite{brenier}. After discretization, this mass transport problem
becomes an assignment problem, which can be solved by efficient
algorithms (see Sec.~4 of Ref.~\cite{mnras} and
Ref.~\cite{bertsekas}). This is the essence of the
Monge--Amp\`ere--Kantorovich (MAK) reconstruction method.

Mohayaee et al.~\cite{royaetal} tested the quality of MAK reconstruction
by applying it to final states of standard N-body simulations, that were
performed for various random initial conditions of cosmological relevance.
The authors of \cite{royaetal} noted the unprecedented accuracy of
the reconstructions down to a few megaparsecs. In particular, they performed
comparisons between three different Lagrangian maps: (i) the map based on
the N-body integration, (ii) the map obtained by applying the MAK procedure to
the current density field, calculated by the N-body integration, (iii) the map
obtained by applying the Zeldovich approximation, starting from the same
initial condition as for the N-body simulation. The conclusion of their
comparisons is that the N-body map is approximated much better by the MAK-generated
map than by the Zeldovich map. This is particularly striking in their Fig.~7,
which shows the negative Lagrangian divergence of the displacement $\x-\q$,
obtained by the three methods mentioned above.

Can we understand this good performance of MAK reconstruction on
sufficiently large spatial scales? First, let us make the rather
obvious observation that any Lagrangian map that is the (Lagrangian)
gradient of a convex function (here called the ``Brenier property'')
will be reconstructed exactly (in a discretized version), if we solve
the associated quadratic-cost optimal transport problem, for example,
by using the MAK procedure. Here is a trivial example of this: if we let an
initially quasi-uniform mass distribution evolve by pure Zeldovich
dynamics to a final distribution and there is no shell crossing, then
the MAK reconstruction is exact. If the solutions to the
Euler--Poisson equations had the Brenier property, MAK would perform
an exact reconstruction, but they don't: flows which
solve the Euler--Poisson equations have no Eulerian vorticity but do
generate Lagrangian vorticity \cite{varenna}.

However there exists a
refinement of the Zeldovich approximation which possesses the Brenier
property. This is given by the second-order of the Lagrangian perturbation theory
\cite{moutarde91,buchert92,buchert94,bouchetal95,catelan95,sahni-shandarin1996}.
Here we shall not describe the Lagrangian perturbation theory in any
technical details since the reader can find them in the above
publications. Nevertheless, in order to discuss some of its conceptual
problems, we describe briefly a few key steps.
One rewrites the Euler--Poisson equations in Lagrangian coordinates, to
obtain a set of nonlinear equations for the displacement $\x-\q$ and
its space and time partial derivatives. Assuming then that in a
suitable sense (see below) the displacement is small, O$(\epsilon)$,
one expands the equations in powers of the small parameter $\epsilon$.
Here, the only perturbed quantities are the deviations of the particle
trajectory from the homogeneous Hubble flow, i.e., from a purely
expanding Universe. At the first order, O$(\epsilon)$, one has the
Zeldovich approximation, which, as discussed in
Sec.~\ref{sub:linea}, is omni-potential. In particular, the Lagrangian
map is potential (by the Brenier property) and the velocity is
potential in both Eulerian and Lagrangian coordinates. As we have
said, with the Lagrangian perturbation theory, one may refine the
Zeldovich approximation. The second-order Lagrangian perturbation theory
(usually denoted by L2) captures significant gravitational physics, for example
some tidal effects \cite{buchertehlers93}, whose importance in the large scale
structure formation has been widely recognized (see, e.g.,
\cite{peeblesgroth76}). L2 has the remarkable property that, for
standard cosmological initial conditions, the Lagrangian map is still
potential. As a consequence, the velocity is also potential in
Lagrangian coordinates. With L2, in Eulerian coordinates the velocity is
potential only up to second order. One would have to sum the whole series
to arbitrarily high orders of the Lagrangian perturbation theory (i.e.,
to arbitrary orders of $\epsilon$) to recover the Eulerian potentiality of the
velocity (but note that convergence of the asymptotic series is not guaranteed
\cite{sahni-shandarin1996,nadkarnigosh-chernoff2011}). L2 having the Brenier property,
the Lagrangian and inverse Lagrangian maps can be reconstructed exactly as
an optimal transport problem, for example, by the MAK technique. This is
probably the main reason why MAK performs well (at sufficiently large scales).
Beyond the second order of the Lagrangian
perturbation theory, scales below the non-linearity scale are
expected to play a decisive role and Lagrangian vorticity is
unavoidable \cite{buchert94}. Such scales cannot be handled accurately
by standard MAK reconstruction.

Finally, let us discuss briefly the thorny issue of the validity of the
Lagrangian perturbation theory. As mentioned before, the
Euler--Poisson equations are a consequence of the Vlasov--Poisson
equations only as long as multi-streaming is absent. The problem is
that, with any cosmologically realistic initial condition at
decoupling, multi-streaming appears immediately or, anyway, well before the
present epoch. The situation is somewhat similar to what we would have
with a one-dimensional Burgers flow in which the initial velocity
would be spatially non-differentiable or nearly so, for example, a
velocity field whose space dependence is the Ornstein--Uhlenbeck
process: shocks would then appear after an arbitrarily short time. This, of
course, constitutes no reason to discard the Euler--Poisson equations for
cosmology. One way to keep them is to \textit{regularize} the initial
velocity field by applying a low-pass filter, i.e., by setting to zero all
Fourier amplitudes whose wavenumbers exceed a cutoff $K =1/L$, where
$L$ is the regularization scale. It may then be shown that no
shell-crossing and thus no multi-streaming will take place for at least
some finite time $T_\star(L)$. Roughly, $T_\star(L)$ is the inverse of
the largest shear (velocity gradient) of the regularized initial
velocity. Up to this time the Euler--Poisson equations are valid.
Furthermore, it may be shown that, up to time $T$, the small parameter which
controls the validity of the Lagrangian perturbation theory (including that
of the Zeldovich approximation) is $T/ T_\star(L)$.
Of course, the true problem is not regularized, or only barely so. Comparison
with N-body simulations suggests nevertheless that the Lagrangian
perturbation theory remains valid for $T/T_\star(L)\ll 1$ when looking only
at scales greater than $L$. Understanding this from a theoretical point of view
is a challenge, which has been partially addressed by Buchert \cite{varenna}.

\section{Concluding remarks}
\label{s:conclusion}

The main question that we have addressed in this paper concerns
omni-potentiality, the (convex) potential character of the mapping
from any time $t_1$ to any time $t_2 \ge t_1$ with $0\le t_1\le t_2\le
T$. First, we have considered a class of flows of ``Zeldovich type'',
comprised of pure Zeldovich/Burgers flows and those obtained from them
by application of arbitrary nonlinear transformations of the time
variable and arbitrary time-dependent scale factors. Such flows are
trivially omni-potential. So are spherically symmetric flows. We then
have investigated (i) the existence of non-trivial omni-potential
flows, (ii) their genericity: can we prescribe the initial velocity
potential in an arbitrary way?

The flows have been characterized by their Lagrangian maps $\q\mapsto\x(\q,t)=
\nabla_{\q}\Phi(\q,t)$ in terms of the scalar potential $\Phi$. As shown in
Sec.~\ref{sub:commu}, omni-potentiality implies that along any particle
trajectory the Hessians ${\cal H}(\Phi(\q,t))$ commute and thus have common
eigendirections. The field of such eigendirections is prescribed as a function
of the initial (Lagrangian) position $\q$, for example, by the eigendirections
of the Hessian of the initial velocity potential $\phi$. The set of
eigendirections of a real symmetric $d\times d$ matrix depends on $d(d-1)/2$
parameters and can be characterized, for example, by $d(d-1)/2$ of the
invariants discussed in the Appendix. As we try to determine a
single scalar function $\Phi$, the situation is rather different in two and
higher dimensions.

When $d=2$, we have a single invariant expressible as a ratio of
suitable combinations of spatial second derivatives of $\Phi$.
Thus, prescribing the field of the invariant values, $g(\q)$, we obtain
a linear second-order PDE, \rf{gen2d}, for $\Phi$.
For a suitable family of fields $g(\q)$, we have found in
Sec.~\ref{sub:ex2d} non-trivial omni-potential flows that are linear
combinations of homogeneous polynomials, thus ensuring the existence of such
non-trivial flows. Using a WKB method, in Sec.~\ref{s:wkb} we have then been able
to construct omni-potential flows, at least locally in space-time, for arbitrary
smooth $g(\q)$. These flows are actually close to Zeldovich flows with straight
trajectories (but are not of ``Zeldovich type''). Extending this
construction globally in space and avoiding the rapid spatial oscillations
inherent to a WKB method constitute interesting open problems.

When $d=3$, omni-potentiality can be expressed in either of two equivalent
ways. One is to demand the commutation of the Hessians of $\Phi$ with those
of a prescribed $\phi$; this gives $d(d-1)/2$ linear homogeneous second-order
PDEs. The other involves working with the invariants, introduced in the
Appendix, which are rational functions of the entries of the
Hessians of $\Phi$; this gives $d(d-1)/2$ nonlinear second-order PDEs. In both
approaches we have one unknown scalar function, $\Phi$, which has to satisfy
more than one equation. Hence, there is an issue of compatibility of these
equations. However, by restricting the potential $\Phi$ to possess a suitable
finite symmetry group, we have obtained a fairly large class of non-trivial
solutions that are even-degree homogeneous polynomials. Whether non-Zeldovich-type
omni-potential flows exist for an arbitrary smooth $\phi$ remains an open
problem. In dimensions $d>3$ the situation is basically the same.

We have shown in Sec.~\ref{sub:app1} that omni-potentiality of a flow is
equivalent to having at each time a velocity field that is potential in
both Eulerian and Lagrangian coordinates. Such double potentiality
was frequently considered in cosmology. It is of particular relevance
when performing reconstruction by convex optimization techniques such
as the Monge--Amp\`ere--Kantorovich (MAK) procedure. Note that the
Euler--Poisson flow is potential in Eulerian coordinates, but not in
general in Lagrangian coordinates. As discussed in Sec.~\ref{s:cosmo},
the approximate Euler--Poisson flow
obtained by the second-order Lagrangian perturbation theory (L2) is
exactly potential in Lagrangian coordinates --- and thus its inverse
Lagrangian map can be obtained exactly by the
Monge--Amp\`ere--Kantorovich (MAK) procedure --- however, it is only
approximately potential in Eulerian coordinates and thus must be
qualified as an approximately omni-potential flow. In particular,
it does not represent an example of an exactly omni-potential three-dimensional
flow for an arbitrary smooth initial velocity.

We finally wish to mention a concrete open problem of cosmological
interest: as mentioned, MAK gives the exact inverse Lagrangian map for
L2 (and this contributes to explaining why MAK works so well when
tested with N-body simulations). However, in L2, between the Eulerian
position $\x$ and its Lagrangian antecedent $\q$, the trajectory is
not given exactly by the Zeldovich approxima\-tion that has a constant
velocity $(\x-\q)/\tau$ (in the coordinates we used in
Sec.~\ref{s:cosmo}). The actual L2 trajectory is slightly curved and
its current (peculiar) velocity differs slightly from
$(\x-\q)/\tau$. It would be of interest to find~how to perturbatively
handle such discrepancies. Given that the full Euler--Poisson
reconstruction problem and the Zeldovich approximation to it both have
convex optimi\-zation formulations, the question arises, whether L2
and higher-order approximations possess such formulations.

\begin{acknowledgments}
We are particularly grateful to J.~Bec, Y.~Brenier, T.~Buchert, S.~Colombi and
A.~Sobolevskii for extensive fruitful discussions. Thanks are also due to
F.~Bouchet, K.~Khanin, R.~Mohayaee, A.~Nusser and E.B.~Vinberg. UF, OP and VZ were
supported by the grant \hbox{ANR-07-BLAN-0235} OTARIE from Agence Nationale de
la Recherche, France. OP and VZ were supported by the grant 11-05-00167-a from
the Russian foundation for basic research. Several visits of BV, OP and VZ to
Observatoire de la C\^ote d'Azur (France) were
supported by the French Ministry of Higher Education and Research.

\end{acknowledgments}

\appendix
\section*{Appendix: Invariants under variation of eigenvalues of symmetric matrices}
\setcounter{equation}{0}
\renewcommand{\theequation}{A\arabic{equation}}

Here we show how the set of commuting symmetric real $d$-dimensional
matrices, having prescribed eigendirections, can be characterized by
a certain number of invariants, which are rational functions of the matrix
entries. The findings here may be of interest beyond the study of
omni-potential flow. To the best of our knowledge, these results are not
available in the published literature. If the reader is aware of any relevant
reference, kindly inform the authors.

By the theorem on codiagonalizability of commuting symmetric matrices,
in an omni-potential flow, Hessians of the potential of the Lagrangian map must
have the same set of eigendirections at all times along any trajectory.
We therefore ask ourselves a more general question: Suppose, the eigenvalues of a symmetric
matrix are being varied while the eigendirections remain fixed. Which quantities
constructed from the entries of the matrix are unaltered under such
variations? We call such quantities {\em invariants}. For a two-dimensional
matrix $H_{ij}$, as we have seen in Sec.~\ref{s:geometry}, the
eigendirections depend on the single invariant $g=(H_{11}-H_{22})/H_{12}$.
When $d<5$, our question can, in principle, be answered by first solving
the characteristic equation whose roots are the eigenvalues and then determining
the eigendirections. However, the expressions for the eigendirections will
then involve radicals, so that obtaining rational expressions for the invariants
is not easy. Anyway, when $d\ge 5$, the non-solvability of the characteristic
equation by radicals renders this strategy useless. We need therefore a more
practical general algebraic approach to this problem.

\subsection{Construction of invariants in dimension $d$}
\label{suba:general}

Let us assume that all eigenvalues $\lambda_i$ of a symmetric $d\times d$
matrix $H$ are distinct, and thus all eigendirections are uniquely defined.
Description of an arbitrary set of $d$ orthogonal eigendirections in $\R^d$
requires $d(d-1)/2$ parameters: $d$ arbitrary directions require $d(d-1)$
parameters, from which one must subtract the number of orthogonality
conditions, $d(d-1)/2$. We expect therefore that a set of $d(d-1)/2$ suitably
chosen invariants uniquely defines the eigendirections.

The problem of finding such invariants is an instance of a much more general
problem of characterizing linear subspaces of a vector space; here the vector
space is that of all real symmetric matrices and the subspace that of matrices
having prescribed eigendirections, which is spanned by the set of all the
powers, from zero to $d-1$, of this matrix. The general problem can be, in
principle, handled using Pl\"ucker coordinates \cite{hodgepedoe}. For our
problem, a more direct approach is available, as now explained. For $d>3$,
our characterization involves fewer invariants than the corresponding number of
Pl\"ucker coordinates.

We denote by ${\bf h}(\lambda_i)$ an eigenvector associated with the eigenvalue
$\lambda_i$, and assume without any loss of generality that no component
of any eigenvector vanishes; one can always achieve this by suitably rotating
the orthonormal basis in $\R^d$, in which the eigenvectors are decomposed.

We construct the invariants as follows: We set,
for some $1\le m\ne n\le d$ and $k\le d$,
\begin{equation}
\gamma^{(d,k)}_{mn}\equiv P^{(d,k)}(\beta_{mn,1},...,\beta_{mn,d}),
\label{hessum}\end{equation}
where
$$\beta_{mn,i}\equiv h_m(\lambda_i)/h_n(\lambda_i)$$
and $P^{(d,k)}$ denote symmetric homogeneous polynomials of degree $k\le d$,
$$P^{(d,k)}(\bm y)\equiv\sum_{1\le j_1<...<j_l<...<j_k\le d}
y_{j_1}...y_{j_l}...y_{j_k}$$
for $\bm y\in\R^d$.
By construction, the quantities $\gamma^{(d,k)}_{mn}$ depend only on
the eigendirections (through the ratios of components) and are invariant under
permutations of the eigendirections; thus they depend only on the set of
eigendirections. Then one substitutes into \rf{hessum} the respective components
of the eigenvectors ${\bf h}(\lambda_i)$, expressed in terms of the associated
eigenvalues $\lambda_i$ and of the entries of the matrix $H$. It is easily
seen that this will produce rational functions of the matrix entries and of
the eigenvalues. Furthermore, it may be shown that the eigenvalues
enter only through symmetric polynomial combinations, which --- by Vi\`ete's
theorem applied to the characteristic polynomial --- have a polynomial
dependence on the matrix entries. The actual derivation of the invariants
can be partially simplified by making use of the identity
\begin{equation}
\prod_{k=1}^d(\lambda_k+c)=\det\|H+cI\|.
\label{prod}\end{equation}

\subsection{Relations between invariants}
\label{suba:relations}

We obtain thus $d^2(d-1)$ invariants $\gamma^{(d,k)}_{mn}$ in the form of
rational functions of the entries of the symmetric matrix $H$. Evidently,
these invariants are too numerous to be independent. For instance, for any $1\le m\ne n\ne l\le d$ and
$0<k<d$ they clearly satisfy the relations (no summation on repeated indices!)
\begin{equation}
\gamma^{(d,d)}_{mn}\gamma^{(d,d)}_{nm}=1,
\label{prodrel}\end{equation}
\begin{equation}
\gamma^{(d,d)}_{ml}\gamma^{(d,d)}_{ln}=\gamma^{(d,d)}_{mn}
\label{prodrel3}\end{equation}
and
\begin{equation}
\gamma^{(d,k)}_{mn}=\gamma^{(d,d)}_{mn}\gamma^{(d,d-k)}_{nm}.
\label{hessdk}\end{equation}
Identities \rf{prodrel} and \rf{hessdk} link invariants $\gamma^{(d,k)}_{mn}$
for, say, $m<n$ with those for $m>n$; there are $d^2(d-1)/2$ of such
independent relations between invariants. Equations \rf{prodrel3} imply
$$\gamma^{(d,d)}_{mn}=\prod_{i=n}^{m-1}\gamma^{(d,d)}_{i+1,i}$$
for any $m>n+1$; conversely, any relation \rf{prodrel3} follows from
these relations together with \rf{prodrel}. Thus, the identities \rf{prodrel3}
contribute further $(d-1)(d-2)/2$ independent relations.

For any $n$ such that $1\le n\le d$, the relation
$$\sum_{\mbox{$\scriptsize m=1\atop m\ne n$}}^d\gamma^{(d,2)}_{mn}=-{d(d-1)\over2}$$
stems from orthogonality of the eigendirections (there are $d$ such relations).

Another family of relations involves the invariants $\gamma^{(d,1)}_{mn}$.
For $p> 0$, let $H^{(p)}$ denote the $p$th power of the matrix $H$ and let
$H^{(p)}_{mn}$ denote its entries, which are of course readily expressed in terms
of the entries of the matrix $H$; we also set $H^{(0)}\equiv I$. The relations
$$\sum_{m=1}^dH^{(p)}_{mn}h_m(\lambda_i)=\lambda_i^ph_n(\lambda_i),$$
that hold true, for any $i$ and $n$, by definition of eigenvectors of $H$,
${\bf h}(\lambda_i)$, and the identity
$$\sum_{i=1}^d\lambda_i^p={\rm tr}H^{(p)}$$
imply, for each $p>0$, the relation
\begin{eqnarray}
&&\sum_{n=1}^d\sum_{\mbox{$\scriptsize m=1\atop m\ne n$}}^dH^{(p)}_{mn}\gamma^{(d,1)}_{mn}\nonumber\\
&=&\sum_{i=1}^d\sum_{n=1}^d\sum_{\mbox{$\scriptsize m=1\atop m\ne n$}}^dH^{(p)}_{mn}
{h_m(\lambda_i)\over h_n(\lambda_i)}\nonumber\\
&=&\sum_{i=1}^d\sum_{n=1}^d(\lambda_i^p-H^{(p)}_{nn})=0.
\label{powH}\end{eqnarray}
Since, by the Cayley--Hamilton theorem, any matrix is a root of its characteristic
polynomial, the entries $H^{(p)}_{mn}$ for $p\ge d$ are linear combinations
of $H^{(p')}_{mn}$ for \hbox{$p-d\le p'\le p-1$}, the coefficients in these linear
combinations being independent of indices $p$, $m$ and $n$. Therefore, relation
\rf{powH} for $p\ge d$ is a consequence of $d$ such relations for $p=p'$ such
that $p-d\le p'\le p-1$. The relation \rf{powH} for $p=0$ is trivial, and hence
there are $d-1$ independent relations \rf{powH} for $1\le p\le d-1$.

In principle, the total number $d^2(d-1)$ of invariants
$\gamma^{(d,k)}_{mn}$ should exceed the number $d(d-1)/2$ of independent
invariants by the number of relations between the invariants.
For arbitrary $d$, we have obtained above
$${d^2(d-1)\over2}+{(d-1)(d-2)\over2}+d+(d-1)={d(d^2+1)\over2}$$
independent relations. For $d=3$, they constitute 15~relations constraining
the 18 invariants that we have introduced; this fits our expectations that
at most 3 invariants are independent, because the orthogonal frame
of eigendirections of a $3\times3$ symmetric matrix is described by 3 Euler
angles. For $d>3$ the number of the obtained independent relations is still
insufficient to fill the gap, and more relations are to be identified.

\subsection{Invariants for $d=2$}
\label{a:R2}

As an example, we present a detailed derivation of the invariant
$\gamma^{(2,1)}_{12}$ for $d=2$. The $i$th eigenvector of a $2\times2$
symmetric matrix $H$ is $(H_{12},\lambda_i-H_{11})$, and hence
$$\gamma^{(2,1)}_{12}\equiv{H_{12}\over\lambda_1-H_{11}}+{H_{12}\over\lambda_2-H_{11}}
={H_{12}\,(\lambda_1+\lambda_2-2H_{11})\over(\lambda_1-H_{11})(\lambda_2-H_{11})}.$$
The characteristic equation for the eigenvalues is
$$\lambda^2-(H_{11}+H_{22})\lambda+H_{11}H_{22}-H_{12}^2=0,$$
and hence $\lambda_1+\lambda_2=H_{11}+H_{22}$. By virtue of \rf{prod},
$$(\lambda_1-H_{11})(\lambda_2-H_{11})=\det\left\|
\begin{array}{cc}
0&H_{12}\\
H_{12}&H_{22}-H_{11}
\end{array}\right\|=-H_{12}^2.$$
Consequently,
$$\gamma^{(2,1)}_{12}={H_{11}-H_{22}\over H_{12}}.$$
The invariant $\gamma^{(2,1)}_{21}$ can be found by a similar calculation, or
just by swapping subscripts in the expression for $\gamma^{(2,1)}_{12}$; clearly,
the two invariants are interrelated: $\gamma^{(2,1)}_{21}=-\gamma^{(2,1)}_{12}$.
The invariant $\gamma^{(2,2)}_{21}$ turns out to be degenerate:
$$\gamma^{(2,2)}_{21}={H^2_{12}\over(\lambda_1-H_{11})
(\lambda_2-H_{11})}=-1.$$
Thus, we have obtained the same invariant as that found in Sec.~\ref{sub:eqn2d}.

\subsection{Invariants for $d=3$}
\label{a:R3}

We consider now invariants in the three-dimensional space.
The $i$th eigenvector of a $3\times3$ symmetric matrix $H$ has components
$$(H_{12}H_{23}+H_{13}(\lambda_i-H_{22}),
\ \ H_{12}H_{13}+H_{23}(\lambda_i-H_{11}),$$
$$(\lambda_i-H_{11})(\lambda_i-H_{22})-H^2_{12}).$$
The procedure outlined above yields
\begin{eqnarray}
\hspace*{-8mm}\gamma^{(3,1)}_{21}\!\!\!\!\!\!&&={H_{22}-H_{11}\over H_{12}}\nonumber\\
&&+\,{H_{13}\over H_{12}}\,{(H_{11}-H_{22})H_{13}H_{23}+(H_{23}^2-H_{13}^2)H_{12}\over
(H_{22}-H_{33})H_{12}H_{13}+(H_{13}^2-H_{12}^2)H_{23}}\nonumber\\
&&+\,{(H_{11}-H_{33})H_{12}H_{23}+(H_{23}^2-H_{12}^2)H_{13}\over
(H_{22}-H_{33})H_{12}H_{13}+(H_{13}^2-H_{12}^2)H_{23}}.
\label{inv3d1}\end{eqnarray}
$\gamma^{(3,3)}_{21}$ is the ratio of two polynomials, which we calculate
using \rf{prod}:
\begin{equation}
\gamma^{(3,3)}_{21}=-{(H_{11}-H_{33})H_{12}H_{23}+
(H_{23}^2-H_{12}^2)H_{13}\over
(H_{22}-H_{33})H_{12}H_{13}+(H_{13}^2-H_{12}^2)H_{23}}.
\label{inv3d3}\end{equation}
$\gamma^{(3,2)}_{21}$ can be found
 by applying identity \rf{hessdk}:
\begin{equation}
\gamma^{(3,2)}_{21}=\gamma^{(3,3)}_{21}\gamma^{(3,1)}_{12}
\label{inv3d2}\end{equation}
(here $\gamma^{(3,1)}_{12}$ can be obtained by permuting the subscripts 1 and 2
in \rf{inv3d1}).

The three invariants $\gamma^{(3,k)}_{21}$ for $1\le k\le 3$ uniquely define
the three ratios $\beta_{21,i}$: by Vi\`ete's theorem, they are roots
of the cubic equation
$$\beta^3-\gamma^{(3,1)}_{21}\beta^2+\gamma^{(3,2)}_{21}\beta^2-\gamma^{(3,3)}_{21}=0.$$
The eigendirections can be recovered in the form of three eigenvectors
$(1,\,\beta_{21,i},\,c_i)$. One can try to obtain the third components $c_i$
($i =1,\,2,\,3)$ from the relations expressing orthogonality
of the eigendirections. However, this produces two solutions:
if $\{c_i\}$ is obtained in this way, then $\{-c_i\}$ is also a solution. Hence,
the invariants $\gamma^{(3,k)}_{21}$, $1\le k\le 3$, define two distinct sets
of eigendirections. The non-uniqueness is eliminated, if in addition we know any
of the invariants $\gamma^{(3,i)}_{j3}$ or $\gamma^{(3,i)}_{3j}$ for $i=1,3$
and $j=1,2$. It is unclear, whether one can choose a set of three invariants
uniquely defining three eigendirections.

When all invariants $\gamma^{(3,k)}_{mn}$ are known for $k=1$ and 3, the equations
for the entries of symmetric matrix $H$ can be considerably simplified.
In view of \rf{inv3d3} and the same equation with permuted subscripts 2 and 3,
relation \rf{inv3d1} can be expressed as
\begin{equation}
\gamma^{(3,1)}_{21}+\gamma^{(3,3)}_{21}={H_{22}-H_{11}\over H_{12}}
+{H_{13}\over H_{12}}\,\gamma^{(3,3)}_{31}.
\label{diag}\end{equation}
Adding \rf{diag} to its analogue, where subscripts 1 and 2 are permuted, we obtain
\begin{eqnarray}
&&-H_{12}\left(\gamma^{(3,1)}_{21}+\gamma^{(3,3)}_{21}+\gamma^{(3,1)}_{12}+\gamma^{(3,3)}_{12}\right)\nonumber\\
&&+\,H_{13}\gamma^{(3,3)}_{31}+H_{23}\gamma^{(3,3)}_{32}=0.
\label{nondiag}\end{eqnarray}
Permuting subscripts in this equation, we obtain a linear system of equations
for the non-diagonal entries of $H$. (Here, we note that the sum of
\rf{nondiag} and its counterparts with permuted subscripts reduces to \rf{powH}
for $p=1$.) Upon solving this linear system, we find the differences $H_{mm}-H_{nn}$ from
\rf{diag} and the analogues of this equation with permuted subscripts, i.e.,
we determine all diagonal entries up to an additive constant. Further
determination of the matrix $H$ would require, of course, the knowledge of its
three eigenvalues. This implies that the non-diagonal entries, as determined
from the above-mentioned linear system, involve two free parameters. This, in
turn, requires that \rf{nondiag} be equivalent to any equation, obtained from
it by permutation of subscripts.

Consequently, we obtain further relations between
the invariants: Permuting subscripts, say, 1 and 3 in \rf{nondiag}, we get
\begin{eqnarray*}
&&-H_{23}\left(\gamma^{(3,1)}_{23}+\gamma^{(3,3)}_{23}+\gamma^{(3,1)}_{32}+\gamma^{(3,3)}_{32}\right)\\
&&+\,H_{13}\gamma^{(3,3)}_{13}+H_{12}\gamma^{(3,3)}_{12}=0.
\end{eqnarray*}
This equation is equivalent to \rf{nondiag} if and only if
$$\left(\gamma^{(3,1)}_{21}+\gamma^{(3,3)}_{21}+\gamma^{(3,1)}_{12}+\gamma^{(3,3)}_{12}\right)
\gamma^{(3,3)}_{13}+\gamma^{(3,3)}_{31}\gamma^{(3,3)}_{12}=0$$
and
\begin{eqnarray*}
&&\hspace*{-15mm}\left(\gamma^{(3,1)}_{21}+\gamma^{(3,3)}_{21}+\gamma^{(3,1)}_{12}+\gamma^{(3,3)}_{12}\right)\\
&&\hspace*{-15mm}\times\left(\gamma^{(3,1)}_{23}+\gamma^{(3,3)}_{23}+\gamma^{(3,1)}_{32}+\gamma^{(3,3)}_{32}\right)
=\gamma^{(3,3)}_{32}\gamma^{(3,3)}_{12}.
\end{eqnarray*}
In view of \rf{prodrel3}, the first of these relations is equivalent to
\begin{equation}
\gamma^{(3,1)}_{21}+\gamma^{(3,3)}_{21}+\gamma^{(3,1)}_{12}+\gamma^{(3,3)}_{12}
+\gamma^{(3,3)}_{31}\gamma^{(3,3)}_{32}=0,
\label{impdouble}\end{equation}
and the second one follows from \rf{impdouble} and the relation obtained from
\rf{impdouble} by permuting subscripts 1 and 3. Relation \rf{impdouble} and its
analogues with permuted subscripts can, of course, be established directly.
Such relations can be used in three dimensions for verifying the consistency
of the invariants, instead of using \rf{powH}, one of our basic relations
between invariants.


\begin{thebibliography}{99}
\bibitem{arnold1972} Arnold, V.I., ``Notes on the behavior of flows of the three-dimensional ideal fluid under a small perturbation of the initial velocity field," Appl. Math. Mech. {\bf 36} 2, 255--262 (1972).
\bibitem{bernardeau} Bernardeau, F., Colombi, S., Gazta{\~n}aga, E., and Scoccimarro, R., ``Large-Scale Structure of the Universe and Cosmological Perturbation Theory," Phys. Rep. {\bf 367}, 1--248 (2002).
\bibitem{bertsekas} Bertsekas, D.P., ``Auction algorithms for network flow problems: A tutorial introduction," Comput. Optim. Appl. {\bf 1}, 7--66 (1992).
\bibitem{bouchetal95} Bouchet, F.R., Colombi, S., Hivon, E., and Juszkiewicz, R., ``Perturbative Lagrangian approach to gravitational instability," Astron. Astrophys. {\bf 296}, 575--608 (1995).
\bibitem{brenierlong} Brenier, Y., ``D\'ecomposition polaire et r\'earrangement
monotone des champs de vecteurs", C. R. Acad. Sci. Paris S\'erie I Math. {\bf 305}, 805--808 (1987).
\bibitem{brenier} Brenier, Y., ``Polar factorization and monotone rearrangement
of vector-valued functions," Comm. Pure Appl. Math. {\bf 44}, 375--417 (1991).
\bibitem{mnras} Brenier,Y., Frisch, U., H\'enon, M., Loeper, G., Matarrese, S., Mohayaee, R., and Sobolevski, A., ``Reconstruction of the early Universe as a convex optimization problem," Mon. Not. R. Astron. Soc. {\bf 346}, 501--524 (2003).
\bibitem{buchert92} Buchert, T., ``Lagrangian theory of gravitational instability of Friedman--Lemaitre cosmologies and the `Zel'dovich approximation'," Mon. Not. R. Astron. Soc. {\bf 254}, 729--737 (1992).
\bibitem{buchert94} Buchert, T., ``Lagrangian theory of gravitational instability of Friedman--Lemaitre cosmologies -- a generic third-order model for non-linear clustering," Mon. Not. R. Astron. Soc. {\bf 267}, 811--820 (1994).
\bibitem{varenna} Buchert, T., ``Lagrangian perturbation approach to the formation of large-scale structure,'' in Proc. IOP Enrico Fermi, Course CXXXII, {\it Dark Matter in the Universe}, Varenna 1995, eds.: S.~Bonometto, J.~Primack, A.~Provenzale, IOS Press Amsterdam, pp. 543--564 (1996).
\bibitem{buchertehlers93} Buchert, T. and Ehlers, J., ``Lagrangian theory of gravitational instability of Friedman--Lemaitre cosmologies-second order approach:an improved model for non-linear clustering," Mon. Not. R. Astron. Soc. {\bf 264}, 375--387 (1993).\\
\bibitem{catelan95} Catelan, P., ``Lagrangian dynamics in non-flat universes
and non-linear gravitational evolution,'' Mon. Not. R. Astron. Soc. {\bf 276}, 115--124 (1995).
\bibitem{nature} Frisch, U., Matarrese, S., Mohayaee, R., and Sobolevski, A.,
``A reconstruction of the initial conditions of the Universe by optimal mass
transportation,'' Nature {\bf 417}, 260--262 (2002).
\bibitem{hodgepedoe} Hodge, W.V.D. and Pedoe, D., {\it Methods of Algebraic
Geometry}, Volume I (Book II). Cambridge University Press (1947).
\bibitem{hornjohnson} Horn, R.A. and Johnson, C.R., {\it Matrix Analysis},
Cambridge University Press, Cambridge (1990).
\bibitem{kantorovich} Kantorovich, L., ``On the translocation of masses,''
C. R. Acad. Sci. URSS {\bf 37}, 199--201 (1942).
\bibitem{loeper} Loeper, G., ``The reconstruction problem for the Euler--Poisson system in cosmology,'' Arch. Rational Mech. Anal., {\bf 179}, 153--216 (2006).
\bibitem{royaetal} Mohayaee, R., Mathis, H., Colombi, S., and Silk, J., ``Reconstruction of primordial density fields,'' Mon. Not. R. Astron.~Soc. {\bf 365}, 939--959 (2006).
\bibitem{monge1781} Monge, G., ``M\'emoire sur la th\'eorie des d\'eblais et
des remblais,'' Hist. Acad. R. Sci. Paris, 666--704 (1781).
\bibitem{moutarde91} Moutarde, F., Alimi, J.M., Bouchet, F.R., Pellat, R., and Ramani, A., ``Precollapse scale invariance in gravitational instability,'' Astrophys. J. {\bf 382}, 377--381 (1991).
\bibitem{nadkarnigosh-chernoff2011} Nadkarni-Ghosh, S. and Chernoff, D.F., ``Extending the domain of validity of the Lagrangian approximation,'' Mon. Not. R. Astron. Soc. {\bf 410}, 1454--1488 (2011).
\bibitem{peebles89} Peebles, P.J.E., ``Tracing Galaxy Orbits Back in Time,'' Astrophys. J. Lett. {\bf 344}, 53--56 (1989).
\bibitem{peeblesgroth76} Peebles, P.J.E. and Groth E.J., ``An integral constraint for the evolution of the galaxy two-point correlation function,'' Astron. Astrophys. {\bf 53}, 131--140 (1976). .
\bibitem{sahni-shandarin1996} Sahni, V. and Shandarin, S., ``Accuracy of Lagrangian approximations in voids,'' Mon. Not. R. Astron. Soc. {\bf 282}, 641--645 (1996).
\bibitem{villani2009} Villani, C., {\it Optimal Transport, Old and New},
Grund\-lehren der mathematischen Wissenschaften, Springer Verlag, Berlin (2009).
\bibitem{zeldovich} Zeldovich (Zel'dovich), Ya.B., ``Gravitational instability: an approximate theory for large density perturbations,'' Astron. Astrophys. {\bf 5}, 84--89 (1970).
\bibitem{sdss} \url{www.sdss.org}
\end{thebibliography}
\end{document}